\documentclass{article}
\newtheorem{Theorem}{Theorem}
\newtheorem{Lemma}[Theorem]{Lemma}

\newtheorem{Definition}[Theorem]{Definition}

\newtheorem{Claim}[Theorem]{Claim}
\newtheorem{Observation}[Theorem]{Observation}
\newcommand{\ex}{\mathrm{ex_{_{\mathcal{P}}}}}
\begin{document}
\title{Extremal $C_{4}$-free/$C_{5}$-free planar graphs}
\author{ Chris Dowden \\
\\
London School of Economics \\
C.Dowden@lse.ac.uk }
\maketitle

\begin{abstract}
We study the topic of `extremal' planar graphs,
defining $\ex(n,H)$ to be the maximum number of edges possible in a planar graph on $n$ vertices
that does not contain a given graph $H$ as a subgraph.
In particular,
we examine the case when $H$ is a small cycle,
obtaining $\ex(n,C_{4}) \leq \frac{15}{7}(n-2)$ for all $n \geq 4$ 
and $\ex(n,C_{5}) \leq \frac{12n-33}{5}$ for all $n \geq 11$,
and showing that both of these bounds are tight. 
\end{abstract}

\section{Introduction}
One of the best-known results in extremal graph theory is Turan's Theorem \cite{tur}, 
which gives the maximum number of edges that a graph on $n$ vertices can have without containing a $K_{r}$ subgraph. 
The Erdos-Stone Theorem \cite{erd} then extends this to the case when $K_{r}$ is replaced by an arbitrary graph $H$, 
showing that the maximum number of edges possible is $(1+o(1)) \left( \frac{\chi (H)-2}{\chi(H)-1} \right) \left( ^{n}_{2} \right)$, 
where $\chi (H)$ denotes the chromatic number of $H$.
This latter result has been called the `fundamental theorem of extremal graph theory' \cite{bol}.

Over the last decade,
a large quantity of work has been carried out in the area of `random' planar graphs
(see, for example, \cite{gim} and~\cite{mcd}).
However,
there seem to be no known results on questions analogous to the Erdos-Stone Theorem,
i.e.~how many edges can a planar graph on $n$ vertices have without containing a given smaller graph?
It is consequently the aim of this paper to now make a start on this topic of `extremal' planar graphs.

Unfortunately,
the case when the forbidden subgraph is a complete graph
(i.e.~the analogue to Turan)
is fairly trivial.
Since $K_{5}$ is not planar,
the only meaningful cases to look at are $K_{3}$ and $K_{4}$,
and these are both straightforward:
for the former,
it can be observed that $K_{2,n-2}$ must be extremal
(since all faces have size four when drawn in the plane),
and so the extremal number of edges is $2n-4$;
for the latter,
it suffices to note that there exist planar triangulations not containing $K_{4}$
(e.g.~take a cycle of length $n-2$
and then add two new vertices that are both adjacent to all those in the cycle),
and so the extremal number is $3n-6$.

The next most natural type of graph to investigate is perhaps a cycle.
Hence,
in this paper we focus on the cases of $C_{4}$ and $C_{5}$,
which turn out to be much more interesting.

We begin,
in Section~\ref{C4},
with the case when the forbidden subgraph is $C_{4}$,
obtaining a tight bound for the extremal number.
In Section~\ref{simpleC5},
we then produce a reasonably simple inequality for $C_{5}$,
before presenting the full $C_{5}$ result in Section~\ref{full},
and demonstrating in Section~\ref{examplessection} that the latter is tight.
In Section~\ref{concluding},
we then finish with some concluding remarks.

\section{$\mathbf{C_{4}}$} \label{C4}

Let us start by introducing a convenient piece of notation:

\begin{Definition}
Let us say that a graph is \emph{H-free}
if it does not contain $H$ as a subgraph
(whether induced or not),
and let \emph{$\ex(n,H)$}
denote the maximum number of edges possible in a planar $H$-free graph on $n$ vertices.
\end{Definition}

In this section,
we shall look at the case when $H=C_{4}$.
In Theorem~\ref{C4thm},
we shall give a simple proof to show that 
$\ex(n,C_{4}) \leq \frac{15}{7}(n-2)$ for all $n \geq 4$,
and then in Theorem~\ref{C4eq}
we shall demonstrate that this bound is actually tight,
in the sense that there are infinitely many values of $n$ for which it is attained exactly.
In the following section,
we shall then move on to looking at $C_{5}$.

We begin with our aforementioned upper bound
(note that in the proof,
as throughout this paper,
we use the term `plane graph'
to mean a particular embedding of a planar graph in the plane):

\begin{Theorem} \label{C4thm}
$\ex(n,C_{4}) \leq \frac{15}{7}(n-2)$ for all $n \geq 4$.
\end{Theorem}
\textbf{Proof}
Let $G$ be an arbitrary $C_{4}$-free connected plane graph on $n \geq 4$ vertices
(clearly it suffices to consider only connected graphs).
Let $f_{i}$ denote the number of faces of size $i$ in $G$,
and let $f$ denote $\sum_{i}f_{i}$.

Note that $\sum_{i \leq 2} f_{i} = f_{4} = 0$
(since $G$ is connected,
$n \geq 4$ and $G \not\supset C_{4}$),
and so 
\begin{eqnarray*}
2e(G) & = & 3f_{3} + \sum_{i \geq 5} if_{i} \\
& \geq & 3f_{3} + 5 \sum_{i \geq 5} f_{i} \\
& = & 3f_{3} + 5(f-f_{3}) \\
& = & 5f - 2f_{3}. 
\end{eqnarray*}
Hence
\begin{equation}
f \leq \frac{2(e(G) + f_{3})}{5}. \label{C4fbound}
\end{equation}

Observe also that no edge of $G$ can be in two faces of size $3$ without creating a $C_{4}$
(unless $G = C_{3}$,
which is not possible if $n \geq 4$),
and so it follows that $e(G) \geq 3f_{3}$.
Putting $f_{3} \leq \frac{e(G)}{3}$ into (\ref{C4fbound}),
we thus obtain
$f \leq \frac{8}{15} e(G)$. 

By Euler's formula,
we then have
$n-2 = e(G) -f \geq \frac{7}{15} e(G)$,
and so $e(G) \leq \frac{15}{7}(n-2)$.
\phantom{qwerty}
\setlength{\unitlength}{0.25cm}
\begin{picture}(1,1)
\put(0,0){\line(1,0){1}}
\put(0,0){\line(0,1){1}}
\put(1,1){\line(-1,0){1}}
\put(1,1){\line(0,-1){1}}
\end{picture} \\

\begin{figure} [ht]
\setlength{\unitlength}{1cm}
\begin{picture}(20,5.45)(-6,-2.4)

\put(0,-0.5){\line(1,2){0.9}}
\put(0,-0.5){\line(-1,2){0.9}}
\put(-1.05,0){\line(1,0){2.1}}
\put(-1.05,0){\line(3,2){1.95}}
\put(1.05,0){\line(-3,2){1.95}}
\put(-0.25,0){\circle*{0.1}}
\put(0.25,0){\circle*{0.1}}
\put(-0.45,0.4){\circle*{0.1}}
\put(0.45,0.4){\circle*{0.1}}
\put(0,0.7){\circle*{0.1}}

\put(-1.05,0){\circle*{0.1}}
\put(1.05,0){\circle*{0.1}}
\put(0,-0.5){\circle*{0.1}}
\put(0.9,1.3){\circle*{0.1}}
\put(-0.9,1.3){\circle*{0.1}}

\put(-1.8,1){\circle*{0.1}}
\put(1.8,1){\circle*{0.1}}
\put(0,2.175){\circle*{0.1}}
\put(1.2,-0.9){\circle*{0.1}}
\put(-1.2,-0.9){\circle*{0.1}}
\put(0,2.2){\line(1,-1){0.9}}
\put(1.8,1){\line(-3,1){0.9}}
\put(1.8,1){\line(-3,-4){0.75}}
\put(1.2,-0.9){\line(-1,6){0.15}}
\put(1.2,-0.9){\line(-3,1){1.2}}
\put(0,2.2){\line(-1,-1){0.9}}
\put(-1.8,1){\line(3,1){0.9}}
\put(-1.8,1){\line(3,-4){0.75}}
\put(-1.2,-0.9){\line(1,6){0.15}}
\put(-1.2,-0.9){\line(3,1){1.2}}

\put(0.4,1.8){\circle*{0.1}}
\put(1.2,1.2){\circle*{0.1}}
\put(1.5,0.6){\circle*{0.1}}
\put(1.1333333,-0.5){\circle*{0.1}}
\put(0.68,-0.7266666){\circle*{0.1}}
\put(-0.4,1.8){\circle*{0.1}}
\put(-1.2,1.2){\circle*{0.1}}
\put(-1.5,0.6){\circle*{0.1}}
\put(-1.1333333,-0.5){\circle*{0.1}}
\put(-0.68,-0.7266666){\circle*{0.1}}

\put(-0.4,1.8){\line(1,0){0.8}}
\put(1.2,1.2){\line(-4,3){0.8}}
\put(1.2,1.2){\line(1,-2){0.3}}
\put(1.1333333,-0.5){\line(1,3){0.3666666}}
\put(1.1333333,-0.5){\line(-2,-1){0.4533333}}
\put(-0.68,-0.7266666){\line(1,0){1.36}}
\put(-1.2,1.2){\line(4,3){0.8}}
\put(-1.2,1.2){\line(-1,-2){0.3}}
\put(-1.1333333,-0.5){\line(-1,3){0.3666666}}
\put(-1.1333333,-0.5){\line(2,-1){0.4533333}}

\put(1.8,3.05){\circle*{0.1}}
\put(3,0){\circle*{0.1}}
\put(0,-2.4){\circle*{0.1}}
\put(-3,0){\circle*{0.1}}
\put(-1.8,3.05){\circle*{0.1}}

\put(0,2.15){\line(2,1){1.8}}
\put(1.8,1){\line(0,1){2.05}}
\put(1.8,1){\line(6,-5){1.2}}
\put(1.2,-0.9){\line(2,1){1.775}}
\put(1.2,-0.9){\line(-4,-5){1.2}}
\put(0,2.15){\line(-2,1){1.8}}
\put(-1.8,1){\line(0,1){2.05}}
\put(-1.8,1){\line(-6,-5){1.2}}
\put(-1.2,-0.9){\line(-2,1){1.8}}
\put(-1.2,-0.9){\line(4,-5){1.2}}

\put(-1.8,3.05){\line(1,0){3.6}}
\put(3,0.05){\line(-2,5){1.2}}
\put(0,-2.4){\line(5,4){3}}
\put(-3,0.05){\line(2,5){1.2}}
\put(0,-2.4){\line(-5,4){3}}

\put(-0.25,0.025){\line(1,0){0.5}}
\put(-0.25,-0.025){\line(1,0){0.5}}
\put(0.25,0.05){\line(1,2){0.2}}
\put(0.25,-0.05){\line(1,2){0.2}}
\put(-0.25,0.05){\line(-1,2){0.2}}
\put(-0.25,-0.05){\line(-1,2){0.2}}
\put(0,0.725){\line(3,-2){0.45}}
\put(0,0.675){\line(3,-2){0.45}}
\put(0,0.725){\line(-3,-2){0.45}}
\put(0,0.675){\line(-3,-2){0.45}}

\put(-0.25,0.035){\line(1,0){0.5}}
\put(-0.25,-0.035){\line(1,0){0.5}}
\put(0.25,0.07){\line(1,2){0.2}}
\put(0.25,-0.07){\line(1,2){0.2}}
\put(-0.25,0.07){\line(-1,2){0.2}}
\put(-0.25,-0.07){\line(-1,2){0.2}}
\put(0,0.745){\line(3,-2){0.45}}
\put(0,0.655){\line(3,-2){0.45}}
\put(0,0.745){\line(-3,-2){0.45}}
\put(0,0.655){\line(-3,-2){0.45}}
\put(0,0.735){\line(3,-2){0.45}}
\put(0,0.665){\line(3,-2){0.45}}
\put(0,0.735){\line(-3,-2){0.45}}
\put(0,0.665){\line(-3,-2){0.45}}

\put(-1.8,3.075){\line(1,0){3.6}}
\put(-1.8,3.025){\line(1,0){3.6}}
\put(3.02,0.05){\line(-2,5){1.22}}
\put(3,0){\line(-2,5){1.2}}
\put(0,-2.425){\line(5,4){3}}
\put(0,-2.375){\line(5,4){3}}
\put(-3.02,0.05){\line(2,5){1.2}}
\put(-3,0){\line(2,5){1.2}}
\put(0,-2.425){\line(-5,4){3}}
\put(0,-2.375){\line(-5,4){3}}

\put(-1.8,3.095){\line(1,0){3.6}}
\put(-1.8,3.005){\line(1,0){3.6}}
\put(3.02,0.035){\line(-2,5){1.22}}
\put(3,-0.03){\line(-2,5){1.2}}
\put(0,-2.445){\line(5,4){3}}
\put(0,-2.355){\line(5,4){3}}
\put(-3.02,0.035){\line(2,5){1.2}}
\put(-3,-0.03){\line(2,5){1.2}}
\put(0,-2.445){\line(-5,4){3}}
\put(0,-2.355){\line(-5,4){3}}
\put(0,-2.435){\line(5,4){3}}
\put(0,-2.365){\line(5,4){3}}
\put(0,-2.435){\line(-5,4){3}}
\put(0,-2.365){\line(-5,4){3}}

\end{picture}

\caption{The icosidodecahedron.} \label{G0}
\end{figure}
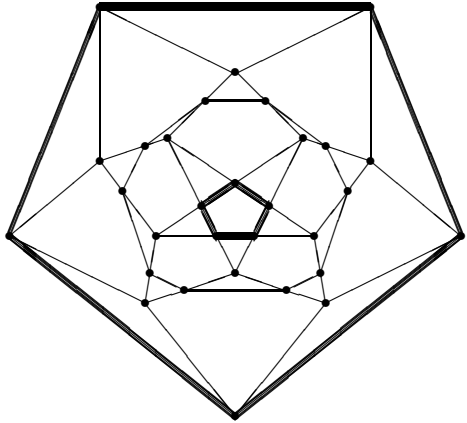

Given the rather uncomplicated nature of the proof of Theorem~\ref{C4thm}
(it only used the facts that $G$ contained no faces of size $4$
and had no edges in two faces of size $3$),
it is perhaps surprising that equality should be achievable.
However,
we shall now see that this is indeed the case:

\begin{Theorem} \label{C4eq}
$\ex(n,C_{4}) = \frac{15}{7} (n-2)$ for $n\equiv 30$ (mod $70$).
\end{Theorem}
\textbf{Proof}
An examination of the proof of Theorem~\ref{C4thm}
shows that equality is achieved for $n$ if and only if
there exists a connected $C_{4}$-free plane graph on $n$ vertices
for which every edge lies in one face of size $3$ and one face of size $5$.
The icosidodecahedron
(see Figure~\ref{G0})
is an example of such a graph with $30$ vertices.

\begin{figure} [ht]
\setlength{\unitlength}{1cm}
\begin{picture}(20,8.75)(-6,-1.75)

\put(0,6){\line(4,-3){4}}
\put(2,0){\line(2,3){2}}
\put(-2,0){\line(1,0){4}}
\put(0,6){\line(-4,-3){4}}
\put(-2,0){\line(-2,3){2}}

\put(-4,3){\line(1,0){8}}
\put(-2,0){\line(1,3){2}}
\put(2,0){\line(-1,3){2}}
\put(-4,3){\line(2,-1){6}}
\put(4,3){\line(-2,-1){6}}

\put(-1.85,0.45){\line(1,0){3.7}}
\put(-1.85,0.45){\line(-2,3){1.7}}
\put(-3.55,3){\line(4,3){3.4}}
\put(-0.15,5.55){\line(1,0){0.3}}
\put(1.85,0.45){\line(2,3){1.7}}
\put(3.55,3){\line(-4,3){3.4}}

\put(-2.15,0.9){\line(1,0){4.3}}
\put(2.933333,2.466666){\line(-4,3){3.597222}}
\put(2.15,0.9){\line(1,2){0.783333}}
\put(-2.933333,2.466666){\line(4,3){3.597222}}
\put(-2.15,0.9){\line(-1,2){0.783333}}

\put(0.663888,5.164583){\line(1,-1){2.164583}}
\put(-0.663888,5.164583){\line(-1,-1){2.164583}}

\put(2.828471,3){\line(-5,-6){2}}
\put(0,0.45){\line(6,1){0.78}}
\put(0.83,0.59){\circle*{0.1}}
\put(-2.828471,3){\line(5,-6){2}}
\put(0,0.45){\line(-6,1){0.78}}
\put(-0.83,0.59){\circle*{0.1}}

\put(0,6){\circle*{0.1}}
\put(2,0){\circle*{0.1}}
\put(4,3){\circle*{0.1}}
\put(-2,0){\circle*{0.1}}
\put(-4,3){\circle*{0.1}}

\put(1.85,0.45){\circle*{0.1}}
\put(0,0.45){\circle*{0.1}}
\put(-1.85,0.45){\circle*{0.1}}

\put(2.15,0.9){\circle*{0.1}}
\put(1.7,0.9){\circle*{0.1}}
\put(0.2,0.9){\circle*{0.1}}
\put(1.078471,0.9){\circle*{0.1}}
\put(-2.15,0.9){\circle*{0.1}}
\put(-1.7,0.9){\circle*{0.1}}
\put(-0.2,0.9){\circle*{0.1}}
\put(-1.078471,0.9){\circle*{0.1}}

\put(0,1){\circle*{0.1}}

\put(1,3){\circle*{0.1}}
\put(2.828471,3){\circle*{0.1}}
\put(3.55,3){\circle*{0.1}}
\put(2.222222,3){\circle*{0.1}}
\put(-1,3){\circle*{0.1}}
\put(-2.828471,3){\circle*{0.1}}
\put(-3.55,3){\circle*{0.1}}
\put(-2.222222,3){\circle*{0.1}}

\put(1.428571,1.714286){\circle*{0.1}}
\put(1.991666,1.995834){\circle*{0.1}}
\put(2.933333,2.466666){\circle*{0.1}}
\put(3.325,2.6625){\circle*{0.1}}
\put(-1.428571,1.714286){\circle*{0.1}}
\put(-1.991666,1.995834){\circle*{0.1}}
\put(-2.933333,2.466666){\circle*{0.1}}
\put(-3.325,2.6625){\circle*{0.1}}

\put(0.592593,4.2222222){\circle*{0.1}}
\put(0.3555555,4.9333333){\circle*{0.1}}
\put(0.15,5.55){\circle*{0.1}}
\put(-0.592593,4.2222222){\circle*{0.1}}
\put(-0.3555555,4.9333333){\circle*{0.1}}
\put(-0.15,5.55){\circle*{0.1}}

\put(0,4.6555555){\circle*{0.1}}
\put(0.663888,5.164583){\circle*{0.1}}
\put(-0.663888,5.164583){\circle*{0.1}}

\put(1.1,0.45){\circle*{0.1}}
\put(-1.1,0.45){\circle*{0.1}}

\put(1.522421,1.43274){\circle*{0.1}}
\put(-1.522421,1.43274){\circle*{0.1}}

\put(2.595299,2.720193){\circle*{0.1}}
\put(-2.595299,2.720193){\circle*{0.1}}

\put(-2.5,7){\line(1,0){5}}
\put(5.5,1){\line(-1,2){3}}
\put(0,-1.75){\line(2,1){5.5}}
\put(2.5,7){\line(-1,0){5}}
\put(-5.5,1){\line(1,2){3}}
\put(0,-1.75){\line(-2,1){5.5}}

\put(-2.5,7.025){\line(1,0){5}}
\put(5.5,1.05){\line(-1,2){3}}
\put(0,-1.775){\line(2,1){5.5}}
\put(2.5,7.025){\line(-1,0){5}}
\put(-5.5,1.05){\line(1,2){3}}
\put(0,-1.775){\line(-2,1){5.5}}

\put(-2.5,6.975){\line(1,0){5}}
\put(5.5,0.95){\line(-1,2){3}}
\put(0,-1.725){\line(2,1){5.5}}
\put(2.5,6.975){\line(-1,0){5}}
\put(-5.5,0.95){\line(1,2){3}}
\put(0,-1.725){\line(-2,1){5.5}}

\put(-2.5,7.045){\line(1,0){5}}
\put(5.5,1.07){\line(-1,2){3}}
\put(0,-1.795){\line(2,1){5.5}}
\put(2.5,7.045){\line(-1,0){5}}
\put(-5.5,1.07){\line(1,2){3}}
\put(0,-1.795){\line(-2,1){5.5}}

\put(-2.5,6.955){\line(1,0){5}}
\put(5.5,0.93){\line(-1,2){3}}
\put(0,-1.705){\line(2,1){5.5}}
\put(2.5,6.955){\line(-1,0){5}}
\put(-5.5,0.93){\line(1,2){3}}
\put(0,-1.705){\line(-2,1){5.5}}

\put(0,-1.785){\line(2,1){5.5}}
\put(0,-1.785){\line(-2,1){5.5}}
\put(0,-1.715){\line(2,1){5.5}}
\put(0,-1.715){\line(-2,1){5.5}}

\put(2.5,7){\circle*{0.1}}
\put(5.5,1){\circle*{0.1}}
\put(0,-1.75){\circle*{0.1}}
\put(-2.5,7){\circle*{0.1}}
\put(-5.5,1){\circle*{0.1}}

\put(-0.3,-0.9){\Large{$G_{0}$}}
\put(-0.5,2){\Large{$G_{k-1}$}}

\put(0,6.025){\line(4,-3){4}}
\put(2,0.0375){\line(2,3){2}}
\put(-2,0.025){\line(1,0){4}}
\put(0,6.025){\line(-4,-3){4}}
\put(-2,0.0375){\line(-2,3){2}}

\put(0,5.975){\line(4,-3){4}}
\put(2,-0.0375){\line(2,3){2}}
\put(-2,-0.025){\line(1,0){4}}
\put(0,5.975){\line(-4,-3){4}}
\put(-2,-0.0375){\line(-2,3){2}}

\put(0,6.045){\line(4,-3){4}}
\put(2,0.0575){\line(2,3){2}}
\put(-2,0.045){\line(1,0){4}}
\put(0,6.045){\line(-4,-3){4}}
\put(-2,0.0575){\line(-2,3){2}}

\put(0,5.955){\line(4,-3){4}}
\put(2,-0.0175){\line(2,3){2}}
\put(-2,-0.005){\line(1,0){4}}
\put(0,5.955){\line(-4,-3){4}}
\put(-2,-0.0175){\line(-2,3){2}}

\end{picture}

\caption{The graph $G_{k}$.} \label{Gk}
\end{figure}
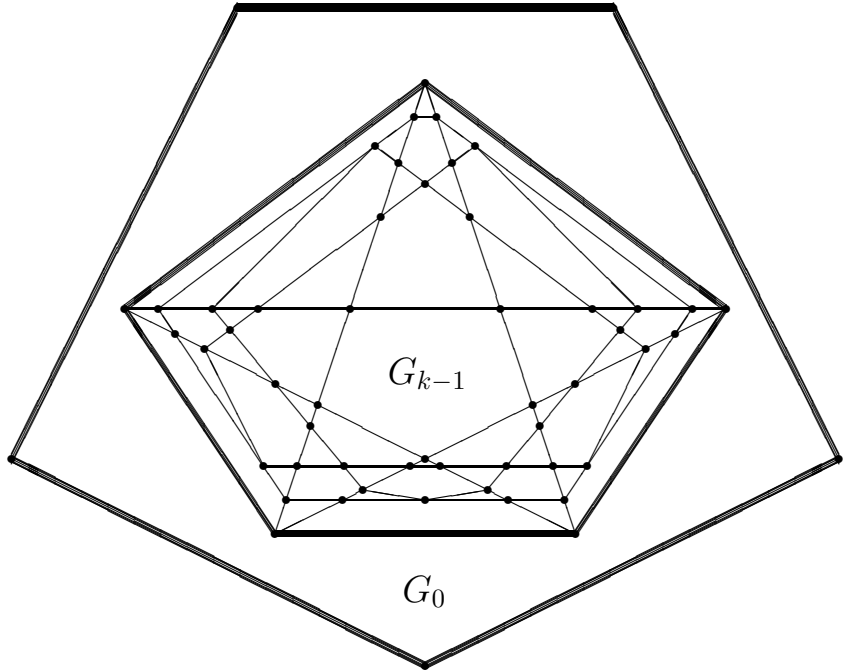

Let us now use the notation $G_{0}$ to denote the icosidodecahedron,
and, for $k \geq 1$,
let us then proceed to define further plane graphs $G_{k}$ of order $30 + 70k$ inductively via the illustration given in Figure~\ref{Gk}.
Here,
the entire graph $G_{k-1}$ is placed into the central pentagon of Figure~\ref{Gk},
and the entire graph $G_{0}$ is then placed between the two bold pentagons of Figure~\ref{Gk}
(in such a way that these are identified with the bold pentagons of Figure~\ref{G0}).

It can be checked that, for all $k$,
$G_{k}$ satisfies the specified conditions,
and so we are done.
\phantom{qwerty}
\setlength{\unitlength}{0.25cm}
\begin{picture}(1,1)
\put(0,0){\line(1,0){1}}
\put(0,0){\line(0,1){1}}
\put(1,1){\line(-1,0){1}}
\put(1,1){\line(0,-1){1}}
\end{picture} \\

As an aside,
we note that slightly more complicated examples can also be produced whenever
$n \equiv 2$ (mod $14$)
(apart from $n \in \{2,16\}$,
for which equality is not possible).

\section{$\mathbf{C_{5}}$ --- a simple inequality} \label{simpleC5}

In Theorem~\ref{C4thm} of the previous section,
we gave a tight bound for $\ex(n,C_{4})$.
In this section,
we shall now move on to looking at $\ex(n,C_{5})$.

By an exactly analogous method to the proof of Theorem~\ref{C4thm}
(with the addition of
taking into account the number of shared edges possible between different $C_{3}$ faces
and between $C_{3}$ faces and $C_{4}$ faces),
we shall first see that
it is reasonably straightforward to obtain $\ex(n,C_{5}) \leq \frac{12}{5}(n-2)$ for all $n \geq 5$
(Lemma~\ref{easier}).
Unfortunately,
we shall then also observe that this time equality is not achievable for large $n$
(Observation~\ref{observation}).
In the following section,
we shall consequently present a much more complicated proof that does provide a tight bound.

Readers may be tempted to skip this section and proceed straight to the full result in Section~\ref{full}.
However, 
it should be noted that the simple inequality given in Lemma~\ref{easier} is actually used during the proof of the full result
(in order to deal with blocks of order $5$ to $10$),
and also that the proof of Observation~\ref{observation} is intended
to aid understanding of the proof of the later result,
which is rather detailed.

We start with the aforementioned `simple' inequality:

\begin{Lemma} \label{easier}
$\ex(n,C_{5}) \leq \frac{12}{5}(n-2)$ for all $n \geq 5$.
\end{Lemma}
\textbf{Proof}
The proof is by induction.
It can be checked that the statement is true for $n=5$,
so let us assume that it is true for all $n \in \{5,6, \ldots, k-1\}$,
and let us now consider an arbitrary $C_{5}$-free plane graph $G$ on $k \geq 6$ vertices.
Let $f_{i}$ denote the number of faces of size $i$ in $G$,
and let $f$ denote $\sum_{i} f_{i}$.
Note that we may assume that $\delta (G) \geq 3$,
since otherwise we could simply delete a vertex of minimal degree 
and obtain the result by applying the induction hypothesis to the resulting graph.

Let us use $e_{3}$ to denote 
the number of edges in $G$
that are in at least one face of size $3$.

\begin{Claim} \label{triclaim}
$f_{3} \leq \frac{e_{3}}{2}$.
\end{Claim}
\textbf{Proof}
Note that each face of size $3$ must be a $C_{3}$.
For each one of these faces,
consider the number of edges that are also in another $C_{3}$ face.
Note that it is possible for there to be two such edges
(if they both belong to the same $K_{4}$ subgraph),
but that it is impossible to have all three
(unless $G=C_{3}$ or $G=K_{4}$,
which may be ignored since $k \geq 6$).
Hence, the total number of edges in two faces of size $3$ is at most $\frac{1}{2}2f_{3}$
(dividing by two to avoid double-counting),
and so we have
\begin{eqnarray*}
f_{3} & = & \frac{1}{3}( \# \textrm{ edges in at least one face of size $3+\#$ edges in two faces of size }3) \\
& \leq & \frac{1}{3} (e_{3} + f_{3}),
\end{eqnarray*}
from which the claim follows. 
\textbf{End of Proof of Claim~\ref{triclaim}} 

\begin{Claim} \label{sqclaim}
$f_{4} \leq \frac{e(G)-e_{3}}{2}$.
\end{Claim}
\textbf{Proof}
Note that each face of size $4$ must be a $C_{4}$
(using the fact that $k>3$),
and recall that each face of size $3$ must be a $C_{3}$.
Observe that it is impossible for a face of size $4$ and a face of size $3$ to have exactly one edge in common
without creating a $C_{5}$,
and that it is impossible for them to have two edges in common without leaving a vertex
in the middle of the facial boundary with degree $2$
(which would contradict the fact that $\delta (G) \geq 3$).
Hence, the faces of size $4$ must be edge-disjoint from the faces of size $3$.
Each edge can be in at most two faces of size $4$,
and each face of size $4$ contains exactly four edges,
so we then obtain $f_{4} \leq \frac{2}{4}(e(G)-e_{3}) = \frac{e(G)-e_{3}}{2}$.
\textbf{End of Proof of Claim~\ref{sqclaim}} \\

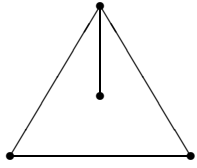
\begin{figure} [ht]
\setlength{\unitlength}{0.8cm}
\begin{picture}(20,2.5)(-6,0)
\put(0,0){\line(1,0){3}}
\put(1.5,2.5){\line(0,-1){1.5}}
\put(0,0){\line(3,5){1.5}}
\put(3,0){\line(-3,5){1.5}}
\put(0,0){\circle*{0.125}}
\put(3,0){\circle*{0.125}}
\put(1.5,1){\circle*{0.125}}
\put(1.5,2.5){\circle*{0.125}}
\end{picture}

\caption{A $C_{3}$ with a pendant edge.} \label{tripen}
\end{figure}

Note that a face of size $5$ in a plane graph must be either a $C_{5}$ or a $C_{3}$ with a pendant edge
(as shown in Figure~\ref{tripen}).
Since $G$ contains no $C_{5}$'s and has minimum degree at least $3$,
we thus have $f_{5}=0$.
Hence (since also $\sum_{i<3}f_{i} = 0$),
we obtain
\begin{eqnarray*}
2e(G) & = & 3f_{3} + 4f_{4} + \sum_{i \geq 6}if_{i} \\
& \geq & 3f_{3} + 4f_{4} + 6(f-f_{3}-f_{4}), 
\end{eqnarray*}
and so
\begin{eqnarray*}
6f & \leq & 2e(G) + 3f_{3} + 2f_{4} \\
& \leq & 2e(G) + \frac{3}{2} e_{3} + e(G) - e_{3} \textrm{ by Claims~\ref{triclaim} and~\ref{sqclaim}} \\
& = & 3e(G) + \frac{1}{2}e_{3} \\
& \leq & \frac{7}{2} e(G), 
\end{eqnarray*}
which gives $f \leq \frac{7}{12}e(G)$.

Using Euler's formula,
we then have $k-2 \geq e(G)-f \geq \frac{5}{12}e(G)$,
from which the result follows.
\phantom{qwerty}
\setlength{\unitlength}{0.25cm}
\begin{picture}(1,1)
\put(0,0){\line(1,0){1}}
\put(0,0){\line(0,1){1}}
\put(1,1){\line(-1,0){1}}
\put(1,1){\line(0,-1){1}}
\end{picture} \\

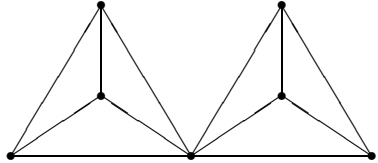
\begin{figure} [ht]
\setlength{\unitlength}{0.8cm}
\begin{picture}(20,2.5)(-4.5,0)
\put(0,0){\line(1,0){3}}
\put(1.5,2.5){\line(0,-1){1.5}}
\put(0,0){\line(3,5){1.5}}
\put(3,0){\line(-3,5){1.5}}
\put(0,0){\circle*{0.125}}
\put(3,0){\circle*{0.125}}
\put(1.5,1){\circle*{0.125}}
\put(1.5,2.5){\circle*{0.125}}
\put(0,0){\line(3,2){1.5}}
\put(3,0){\line(-3,2){1.5}}

\put(3,0){\line(1,0){3}}
\put(4.5,2.5){\line(0,-1){1.5}}
\put(3,0){\line(3,5){1.5}}
\put(6,0){\line(-3,5){1.5}}
\put(6,0){\circle*{0.125}}
\put(4.5,1){\circle*{0.125}}
\put(4.5,2.5){\circle*{0.125}}
\put(3,0){\line(3,2){1.5}}
\put(6,0){\line(-3,2){1.5}}
\end{picture}

\caption{A graph achieving equality in Lemma~\ref{easier} for $n=7$.} \label{n7eq}
\end{figure}

Figure~\ref{n7eq} gives an example of a graph for which equality in Lemma~\ref{easier} is attained with $n=7$.
However,
we shall now see (in Observation~\ref{observation})
that equality is not possible for larger values of $n$.

Readers may proceed straight to Section~\ref{full} if they so wish.
However,
it is hoped that the ideas presented in the proof of Observation~\ref{observation},
and in particular the introduction of a graph $G^{\prime}$,
will prove helpful in understanding the later detailed arguments.

\begin{figure} [ht]
\setlength{\unitlength}{0.8cm}
\begin{picture}(20,2.5)(-4.5,0)
\put(0,0){\line(1,0){3}}
\put(0,0){\line(3,5){1.5}}
\put(3,0){\line(-3,5){1.5}}
\put(0,0){\circle*{0.125}}
\put(3,0){\circle*{0.125}}
\put(1.5,2.5){\circle*{0.125}}

\put(3,0){\line(1,0){3}}
\put(3,0){\line(3,5){1.5}}
\put(6,0){\line(-3,5){1.5}}
\put(6,0){\circle*{0.125}}
\put(4.5,2.5){\circle*{0.125}}
\end{picture}

\caption{The graph $G^{\prime}$ formed from the graph $G$ in Figure~\ref{n7eq}.} \label{obsfig}
\end{figure}
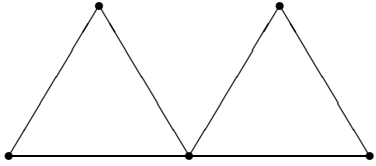

\begin{Observation} \label{observation}
$\ex(n,C_{5}) < \frac{12}{5}(n-2)$ for all $n \geq 8$.
\end{Observation}
\textbf{Proof}
A careful examination of the proof of Lemma~\ref{easier} shows that equality is possible only if 
the plane graph $G$ is connected and consists entirely of $K_{4}$'s and faces of size $6$,
combined in such a way that no two $K_{4}$'s share an edge 
and no two faces of size $6$ share an edge
(note that the outside face in Figure~\ref{n7eq} has size $6$).
It will consequently be convenient for us to consider the plane graph $G^{\prime}$ 
formed from $G$ by deleting the central vertex from each $K_{4}$,
along with all incident edges,
to leave a $C_{3}$
(see Figure~\ref{obsfig} for an example).
Hence, 
$G^{\prime}$ will then consist entirely of faces of size $3$ and $6$,
\emph{with every edge of $G^{\prime}$
lying on the boundary between a face of size $3$ and a face of size $6$}.

Let us use $f^{\prime}$ to denote the total number of faces in $G^{\prime}$,
$f_{3}^{\prime}$ to denote the number of faces in $G^{\prime}$ that have size $3$,
and $f_{6}^{\prime}$ to denote the number of faces that have size $6$.
Then we consequently obtain $3f_{3}^{\prime} = e(G^{\prime}) = 6f_{6}^{\prime}$,
and $f^{\prime} = f_{3}^{\prime} + f_{6}^{\prime} = \frac{e(G^{\prime})}{3} + \frac{e(G^{\prime})}{6} = \frac{e(G^{\prime})}{2}$.
Thus, Euler's formula gives $|G^{\prime}|-2 = e(G^{\prime})-f^{\prime} =  \frac{e(G^{\prime})}{2}$,
and hence $e(G^{\prime}) = 2|G^{\prime}|-4$.

Crucially,
this implies that the average degree in $G^{\prime}$ is less than $4$.
However,
vertices of degree $0$ or $1$ are clearly impossible,
vertices of degree $3$ are also not possible
(every vertex must have even degree,
since the adjacent faces are \emph{alternately} size $3$ and size $6$),
and it can be checked that
a vertex of degree $2$ is only possible for the special case
shown in Figure~\ref{obsfig}
(using the fact that if the faces adjacent to such a vertex $v$  are cycles
$u_{1}vu_{2}$ and $u_{1}vu_{2}u_{3}u_{4}u_{5}$,
then $u_{1}u_{2}u_{3}u_{4}u_{5}$ would be a $C_{5}$).
Hence,
we obtain a contradiction.
\phantom{qwerty}
\setlength{\unitlength}{0.25cm}
\begin{picture}(1,1)
\put(0,0){\line(1,0){1}}
\put(0,0){\line(0,1){1}}
\put(1,1){\line(-1,0){1}}
\put(1,1){\line(0,-1){1}}
\end{picture}

\section{$\mathbf{C_{5}}$ --- the full result} \label{full}

In this section,
we shall improve on the results of Lemma~\ref{easier} and Observation~\ref{observation}
by showing that 
$\ex(n,C_{5}) \leq \frac{12n-33}{5}$ for all $n \geq 11$ (see Theorem~\ref{mainthm}).
As the details are very lengthy,
a sketch of the proof is also provided.
In the following section,
we shall then demonstrate that this bound is tight,
in the sense that there are infinitely many values of $n$
for which it is attained exactly.

\begin{Theorem} \label{mainthm}
$\ex(n,C_{5}) \leq \frac{12n-33}{5}$ for all $n \geq 11$.
\end{Theorem}
\textbf{Sketch of Proof}
Part I:
We take a $C_{5}$-free plane graph $G$,
and use induction to deal with the cases when $\delta(G) \leq 2$ or $\kappa(G) \leq 1$,
where $\kappa(G)$ denotes the vertex-connectivity of $G$
(Lemma~\ref{easier} is utilised for this latter case).
This allows us to simplify some of the details in later arguments. 

Part II:
Similarly to with the proof of Observation~\ref{observation},
we now form a new plane graph $G^{\prime}$
by deleting any edges that lie between two $C_{3}$ faces,
and we again show that $e(G^{\prime}) \leq 2|G^{\prime}|-4$
(this time,
the proof is complicated by the fact that we do not assume that $G$ achieves equality in Lemma~\ref{easier},
and so we need to allow for the presence of faces of size $4$),
and hence that $G^{\prime}$ contains vertices of degree less than $4$.

Part III:
We examine the case when $e(G^{\prime}) = 2|G^{\prime}|-4$,
and deduce that equality implies that either
(i)~$G^{\prime}$ consists solely of $C_{4}$ faces or
(ii)~every edge in $G^{\prime}$ lies on the boundary between one $C_{3}$ face and one $C_{6}$ face.
We then show that the former would imply that $G=G^{\prime}$,
and so $e(G)$ would also only be $2|G|-4$,
while the latter is impossible
(as in the proof of Observation~\ref{observation}). 

Part IV:
We then examine the remaining case when $e(G^{\prime}) \leq 2|G^{\prime}|-5$,
noting that this implies a more considerable lower bound on the number of vertices of degree $2$ or $3$ in $G^{\prime}$.
We find that the presence of each such vertex forces fewer faces of small size 
(to avoid a $C_{5}$),
and we consequently obtain detailed inequalities
which essentially express the exact amount of negative impact that each of these vertices has on $e(G)$.
Crucially,
when combining this with our lower bound for the number of such vertices,
we find that the best we can do is to obtain $e(G) = \frac{12|G|-33}{5}$. \\
\\
\textbf{Full Proof} \\
\textbf{Part I}
Our proof is by induction.
It can be checked systematically that the result $\ex(n,C_{5}) \leq \frac{12n-33}{5}$ holds for $n \in \{11,12,13\}$
(this turns out to be not as onerous as it might at first appear!),
so let $k \geq 14$ and let us assume that the result holds for all $n \in \{11,12, \ldots, k-3,k-2,k-1\}$.

Let $G$ be a $C_{5}$-free plane graph on $k$ vertices,
and let us deal in this part of the proof with the cases when $\delta (G) \leq 2$ or $\kappa (G) \leq 1$.
For the former,
we may simply delete a vertex of minimum degree 
and then use the induction hypothesis to obtain the result.
For the latter,
by considering an arbitrary block in $G$
(which will have $k-i$ vertices for some $i \geq 1$)
and the union of all other blocks,
we may obtain
\begin{eqnarray*}
e(G) & \leq & \max_{i \geq 1} \{ \ex(i+1,C_{5}) + \ex(k-i, C_{5}) \} \\
& \leq & \max_{i \geq 4} \left\{ 1 + \frac{12(k-1)-33}{5}, 3 + \frac{12(k-2)-33}{5}, 6 + \frac{12(k-3)-33}{5}, \right. \\
& & \left. \phantom{wwwl} \frac{12}{5}(i-1) + \frac{12}{5}(k-i-2) \right\} \\
& & \textrm{ using the induction hypothesis and Lemma~\ref{easier}} \\
& = & \frac{12k-36}{5}. 
\end{eqnarray*} 

\begin{figure} [ht]
\setlength{\unitlength}{1cm}
\begin{picture}(20,4) (-1.5,-2.75)
\put(0,0){\line(1,0){1.5}}
\put(0.75,1.25){\line(0,-1){0.75}}
\put(0,0){\line(3,5){0.75}}
\put(1.5,0){\line(-3,5){0.75}}
\put(0,0){\circle*{0.1}}
\put(1.5,0){\circle*{0.1}}
\put(0.75,0.5){\circle*{0.1}}
\put(0.75,1.25){\circle*{0.1}}
\put(0,0){\line(3,2){0.75}}
\put(1.5,0){\line(-3,2){0.75}}

\put(0,0){\line(-1,-2){0.5}}
\put(0,-2){\line(-1,2){0.5}}
\put(1.5,0){\line(1,-2){0.5}}
\put(1.5,-2){\line(1,2){0.5}}
\put(0,-2){\line(1,0){1.5}}

\put(2,-1){\line(1,0){0.75}}
\put(-0.5,-1){\line(-1,0){0.75}}
\put(1.5,0){\line(5,-4){1.25}}
\put(1.5,-2){\line(5,4){1.25}}
\put(0,0){\line(-5,-4){1.25}}
\put(0,-2){\line(-5,4){1.25}}

\put(-1.25,-1){\circle*{0.1}}
\put(-0.5,-1){\circle*{0.1}}
\put(2,-1){\circle*{0.1}}
\put(2.75,-1){\circle*{0.1}}
\put(0,-2){\circle*{0.1}}
\put(1.5,-2){\circle*{0.1}}

\put(0.55,-2.75){\Large{$G$}}

\put(3.75,-1){\vector(1,0){1.5}}

\put(7.5,0){\line(1,0){1.5}}
\put(7.5,0){\line(3,5){0.75}}
\put(9,0){\line(-3,5){0.75}}
\put(7.5,0){\circle*{0.1}}
\put(9,0){\circle*{0.1}}
\put(8.25,1.25){\circle*{0.1}}

\put(7.5,0){\line(-1,-2){0.5}}
\put(7.5,-2){\line(-1,2){0.5}}
\put(9,0){\line(1,-2){0.5}}
\put(9,-2){\line(1,2){0.5}}
\put(7.5,-2){\line(1,0){1.5}}

\put(9,0){\line(5,-4){1.25}}
\put(9,-2){\line(5,4){1.25}}
\put(7.5,0){\line(-5,-4){1.25}}
\put(7.5,-2){\line(-5,4){1.25}}

\put(6.25,-1){\circle*{0.1}}
\put(7,-1){\circle*{0.1}}
\put(9.5,-1){\circle*{0.1}}
\put(10.25,-1){\circle*{0.1}}
\put(7.5,-2){\circle*{0.1}}
\put(9,-2){\circle*{0.1}}

\put(8.05,-2.75){\Large{$G^{\prime}$}}
\end{picture}

\caption{Constructing the graph $G^{\prime}$ from the graph $G$.} \label{Gprime}
\end{figure}
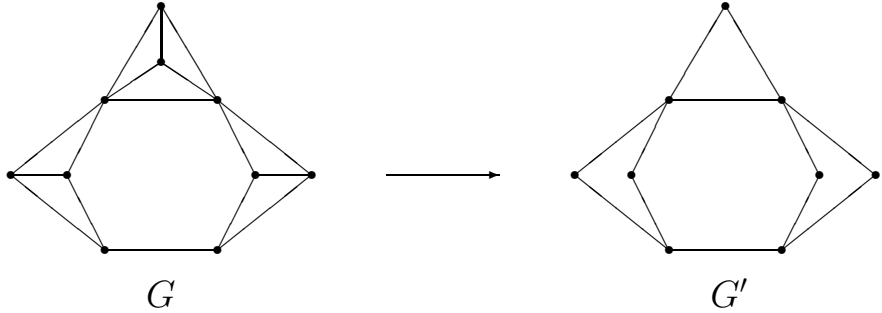

Hence,
for the remainder of the proof,
we may assume that we have
$\kappa (G) \geq 2$
and $\delta (G) \geq 3$.
Note that the condition that there are no cut-vertices implies that all faces of $G$ must be cycles.
We shall often use this fact implicitly in the rest of the proof. \\
\\
\textbf{Part II}
Let $G^{\prime}$ denote the plane graph formed from $G$ by
(a) deleting any edges that lie on two $C_{3}$ faces in $G$
and (b) also deleting any isolated vertices that are thus produced
(see Figure~\ref{Gprime} for an example of this).
In this part of the proof,
we shall work towards showing that
$e(G^{\prime}) \leq 2|G^{\prime}|-4$
by paying particular attention to the number of faces in $G^{\prime}$ that are $C_{3}$'s or $C_{4}$'s.

Observe that if two $C_{3}$ faces in $G$ share an edge,
then (since $G$ is $C_{5}$-free)
these two faces must either be edge-disjoint to all other $C_{3}$ faces in $G$
or must form part of a $K_{4}$ consisting of three adjacent $C_{3}$ faces in $G$.
In the former case,
the two $C_{3}$ faces in $G$ will be replaced by a $C_{4}$ face in $G^{\prime}$;
in the latter case,
the $K_{4}$ in $G$ will be replaced by a $C_{3}$ face in $G^{\prime}$.

Thus, all faces of $G^{\prime}$ will be cycles, as with $G$.
For the rest of the proof,
it will be extremely important to note further that all faces of size greater than $4$ in $G^{\prime}$ were also faces in $G$;
that all $C_{4}$ faces in $G^{\prime}$ were either also faces in $G$ 
or are the remains of two adjacent $C_{3}$ faces of $G$;
and that all $C_{3}$ faces in $G^{\prime}$ were either faces in $G$
or formed the outside of a $K_{4}$ in $G$.
It follows from the latter observation that there will now be no edges that lie on two $C_{3}$ faces
(we may discount the possibility that $G^{\prime} = C_{3}$,
since this would imply that $G=K_{4}$ and hence that $k=4$).

Let $f_{i}^{\prime}$ denote the number of $C_{i}$ faces in $G^{\prime}$,
and let $e_{i}^{\prime}$ denote the number of edges of $G^{\prime}$ that are in at least one $C_{i}$ face.
Then, since there is no double-counting, we have
\begin{equation}
f_{3}^{\prime} = \frac{1}{3} e_{3}^{\prime}. \label{f3prime}
\end{equation}

\begin{figure} [ht]
\setlength{\unitlength}{1cm}
\begin{picture}(20,2.4) (-5.5,-1.1)
\put(1,0){\line(-1,1){1}}
\put(1,0){\line(-1,-1){1}}
\put(-1,0){\line(1,1){1}}
\put(-1,0){\line(1,-1){1}}
\put(0,0){\oval(4,2)[r]}
\put(0,-1){\circle*{0.1}}
\put(0,1){\circle*{0.1}}
\put(-1,0){\circle*{0.1}}
\put(1,0){\circle*{0.1}}
\put(-1.5,-0.05){$v_{1}$}
\put(-0.5,1){$v_{2}$}
\put(-0.5,-1.1){$v_{4}$}
\put(0.5,-0.05){$v_{3}$}
\put(-0.3,-0.1){\large{$F_{1}$}}
\put(1.3,-0.1){\large{$F_{2}$}}
\end{picture}

\caption{A $C_{4}$ face and a $C_{3}$ face that share two edges.} \label{triquadrilateral}
\end{figure}
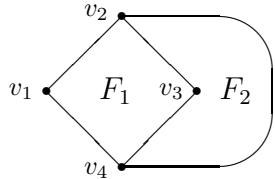

\begin{figure} [ht]
\setlength{\unitlength}{1cm}
\begin{picture}(20,2.4) (-5.5,-1.1)
\put(1,0){\line(-1,1){1}}
\put(1,0){\line(-1,-1){1}}
\put(-1,0){\line(1,1){1}}
\put(-1,0){\line(1,-1){1}}
\put(0,0){\oval(4,2)[r]}
\put(0,-1){\circle*{0.1}}
\put(0,1){\circle*{0.1}}
\put(-1,0){\circle*{0.1}}
\put(1,0){\circle*{0.1}}
\put(-1.5,-0.05){$v_{1}$}
\put(-0.5,1){$v_{2}$}
\put(-0.5,-1.1){$v_{4}$}
\put(0.5,-0.05){$v_{3}$}
\put(1.5,0){\circle*{0.1}}
\put(1,0){\line(1,0){0.5}}
\put(1.5,0){\line(-3,2){1.5}}
\put(1.5,0){\line(-3,-2){1.5}}

\put(1,0.025){\line(-1,-1){1}}
\put(-1,0.025){\line(1,1){1}}
\put(-1,0.025){\line(1,-1){1}}
\put(1,0.025){\line(1,0){0.5}}
\put(1.5,0.025){\line(-3,2){1.5}}

\put(1,-0.025){\line(-1,-1){1}}
\put(-1,-0.025){\line(1,1){1}}
\put(-1,-0.025){\line(1,-1){1}}
\put(1,-0.025){\line(1,0){0.5}}
\put(1.5,-0.025){\line(-3,2){1.5}}

\put(1,0.045){\line(-1,-1){1}}
\put(-1,0.045){\line(1,1){1}}
\put(-1,0.045){\line(1,-1){1}}
\put(1,0.045){\line(1,0){0.5}}
\put(1.5,0.045){\line(-3,2){1.5}}

\put(1,-0.045){\line(-1,-1){1}}
\put(-1,-0.045){\line(1,1){1}}
\put(-1,-0.045){\line(1,-1){1}}
\put(1,-0.045){\line(1,0){0.5}}
\put(1.5,-0.045){\line(-3,2){1.5}}
\end{picture}

\caption{A $C_{5} \subset G$ in case (a).} \label{caseaC5}
\end{figure}
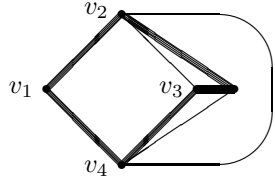

\begin{Claim} \label{triquadrilateralclaim}
A $C_{4}$ face in $G^{\prime}$ cannot share an edge with a $C_{3}$ face.
\end{Claim}
\textbf{Proof}
Clearly, a $C_{4}$ face cannot share exactly one edge with a $C_{3}$ face without creating a $C_{5}$.
If a $C_{4}$ face shares two edges with a $C_{3}$ face,
then it must look as shown in Figure~\ref{triquadrilateral},
where a $C_{4}$ face $F_{1} = v_{1}v_{2}v_{3}v_{4}$
and a $C_{3}$ face $F_{2} = v_{2}v_{3}v_{4}$ share the edges $v_{2}v_{3}$ and $v_{3}v_{4}$.
However, note that if such a situation were to occur in $G^{\prime}$, 
then (since $\deg_{G}(v_{3}) \geq 3$) it must be that either
(a) $F_{2}$ formed the outside of a $K_{4}$ in $G$,
in which case we could then find a $C_{5} \subset G$ (see Figure~\ref{caseaC5}),
or (b) $F_{1}$ is the remains of two $C_{3}$ faces $v_{1}v_{2}v_{3}$ and $v_{1}v_{3}v_{4}$,
in which case the vertex $v_{3}$ and all edges incident to it should have been deleted when $G^{\prime}$ was constructed.
Thus, it must be that such a situation does not occur.
\textbf{End of Proof of Claim~\ref{triquadrilateralclaim}} \\

Note that $f_{4}^{\prime} \leq \frac{e_{4}^{\prime}}{2}$,
with equality if and only if every edge in a $C_{4}$ face is actually in two such faces.
From Claim~\ref{triquadrilateralclaim},
we know that all $C_{4}$ faces in $G^{\prime}$ must be edge-disjoint from all $C_{3}$ faces,
and so it follows that
\begin{equation}
f_{4}^{\prime} \leq \frac{e(G^{\prime}) - e_{3}^{\prime}}{2}, \label{f4prime}
\end{equation}
with equality if and only if every edge not in a $C_{3}$ face is in two $C_{4}$ faces.

Let $f^{\prime}$ denote $\sum_{i}f_{i}^{\prime}$.
Then we have
\begin{eqnarray*}
2e(G^{\prime}) & \geq & 3f_{3}^{\prime} + 4f_{4}^{\prime} + 6(f^{\prime} - f_{3}^{\prime} - f_{4}^{\prime})
\end{eqnarray*}
and so
\begin{eqnarray*}
6f^{\prime} & \leq & 2e(G^{\prime}) + 3f_{3}^{\prime} + 2f_{4}^{\prime} \\
& \leq & 2e(G^{\prime}) + e_{3}^{\prime} + e(G^{\prime}) - e_{3}^{\prime} \textrm{ by (\ref{f3prime}) and (\ref{f4prime})} \\
& = & 3e(G^{\prime}).
\end{eqnarray*}
Thus, $f^{\prime} \leq \frac{1}{2}e(G^{\prime})$.

Hence, since Euler's formula then gives $e(G^{\prime}) = |G^{\prime}| - 2 + f^{\prime} \leq |G^{\prime}| - 2 + \frac{1}{2}e(G^{\prime})$,
we obtain $\frac{1}{2}e(G^{\prime}) \leq |G^{\prime}| - 2$,
i.e.~$e(G^{\prime}) \leq 2|G^{\prime}| - 4$. \\
\\
\textbf{Part III} 
In this part of the proof,
we shall deal with the case when
$e(G^{\prime}) = 2|G^{\prime}| - 4$.
As mentioned in the sketch of the proof,
we shall find that our argument splits into two subcases depending on whether 
$G^{\prime}$ consists solely of $C_{4}$ faces
or is instead comprised of a mixture of $C_{3}$'s and $C_{6}$'s.

Let us start by noting
(from an examination of Part II)
that equality in $e(G^{\prime}) \leq 2|G^{\prime}|-4$
implies that $f^{\prime} = f_{3}^{\prime} + f_{4}^{\prime} + f_{6}^{\prime}$
and that every edge of $G^{\prime}$ not in a $C_{3}$ face is in two $C_{4}$ faces.

Let $E_{i,j}^{\prime}$ denote the set of edges of $G^{\prime}$ between a $C_{i}$ face and a $C_{j}$ face,
and let $E_{i,i}^{\prime}$ denote the set of edges of $G^{\prime}$ between two $C_{i}$ faces.
Recall that $E_{3,3}^{\prime} = \emptyset$,
and also that $E_{3,4}^{\prime} = \emptyset$.
Hence, $E_{3,6}^{\prime} \cup E_{4,4}^{\prime} = E(G^{\prime})$.

If $E_{3,6}^{\prime}$ and $E_{4,4}^{\prime}$ are both non-empty,
then there exists a face with edges in both $E_{3,6}^{\prime}$ and $E_{4,4}^{\prime}$, which is clearly absurd.
Thus, either $E_{3,6}^{\prime} = E(G^{\prime})$ or $E_{4,4}^{\prime} = E(G^{\prime})$.

Let us first consider the case $E_{4,4}^{\prime} = E(G^{\prime})$,
when every edge in $G^{\prime}$ is shared by two $C_{4}$ faces.

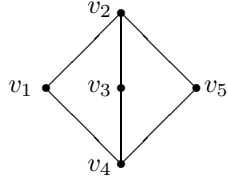
\begin{figure} [ht]
\setlength{\unitlength}{1cm}
\begin{picture}(20,2.4) (-6,-1.1)
\put(1,0){\line(-1,1){1}}
\put(1,0){\line(-1,-1){1}}
\put(-1,0){\line(1,1){1}}
\put(-1,0){\line(1,-1){1}}
\put(0,-1){\line(0,1){2}}
\put(0,-1){\circle*{0.1}}
\put(0,1){\circle*{0.1}}
\put(-1,0){\circle*{0.1}}
\put(1,0){\circle*{0.1}}
\put(0,0){\circle*{0.1}}
\put(-1.5,-0.05){$v_{1}$}
\put(-0.45,1){$v_{2}$}
\put(-0.45,-1.1){$v_{4}$}
\put(1.1,-0.05){$v_{5}$}
\put(-0.45,-0.05){$v_{3}$}
\end{picture}
\caption{Two $C_{4}$ faces that share two edges.} \label{quadrilateraltwo}
\end{figure}

\begin{Claim} \label{quadrilateralquadrilateral}
A $C_{4}$ face in $G^{\prime}$ can only share an edge with another $C_{4}$ face 
if neither is the remains of two adjacent $C_{3}$ faces of $G$.
\end{Claim}
\textbf{Proof}
First, note that we may discount the possibility that $G^{\prime} = C_{4}$,
since this would imply that $G$ also had only $4$ vertices.
Hence, two $C_{4}$ faces in $G^{\prime}$ can only share at most two edges.
In the case when two $C_{4}$ faces share exactly one edge,
it is clear that neither can be the remains of two adjacent $C_{3}$ faces without creating a $C_{5}$ in $G$.
If two $C_{4}$ faces share two edges,
then they must look as shown in Figure~\ref{quadrilateraltwo},
where the $C_{4}$ faces $v_{1}v_{2}v_{3}v_{4}$ and $v_{2}v_{3}v_{4}v_{5}$ share the edges $v_{2}v_{3}$ and $v_{3}v_{4}$.
However, note that the condition $\delta (G) \geq 3$ implies that we must then have either $v_{1}v_{3} \in E(G)$ or $v_{3}v_{5} \in E(G)$,
and both of these would produce a $C_{5}$.
\textbf{End of Proof of Claim~\ref{quadrilateralquadrilateral}} \\

From Claim~\ref{quadrilateralquadrilateral},
it follows that if $E_{4,4}^{\prime} = E(G^{\prime})$ then there can be no edges in $E(G^{\prime}) \setminus E(G)$,
and hence $E(G) = 2k-4 < \frac{12k - 33}{5}$,
and so we would be done.

\begin{figure} [ht]
\setlength{\unitlength}{1cm}
\begin{picture}(20,2.8) (-5.25,-2.3)
\put(0,0){\line(1,0){1.5}}
\put(0,0){\circle*{0.1}}
\put(1.5,0){\circle*{0.1}}
\put(0,0){\line(-1,-2){0.5}}
\put(0,-2){\line(-1,2){0.5}}
\put(1.5,0){\line(1,-2){0.5}}
\put(1.5,-2){\line(1,2){0.5}}
\put(0,-2){\line(1,0){1.5}}
\put(-0.5,-1){\circle*{0.1}}
\put(2,-1){\circle*{0.1}}
\put(0,-2){\circle*{0.1}}
\put(1.5,-2){\circle*{0.1}}
\put(0,-1){\oval(3,2)[l]}
\put(-0.2,0.2){$v_{1}$}
\put(1.3,0.2){$v_{3}$}
\put(-0.85,-1.05){$v$}
\put(2.1,-1.05){$v_{4}$}
\put(-0.2,-2.3){$v_{2}$}
\put(1.3,-2.3){$v_{5}$}
\end{picture}

\caption{A vertex $v$ of degree 2 in the case when $E_{3,6}^{\prime} = E(G^{\prime})$.} \label{trihexdeg2}
\end{figure}
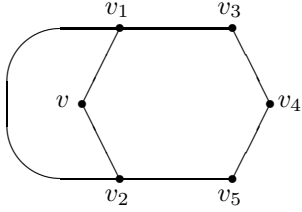

Let us now consider the case when $E_{3,6}^{\prime} = E(G^{\prime})$,
i.e.~when all edges in $G^{\prime}$ lie on the boundary between a $C_{3}$ face and a $C_{6}$ face.
In particular, let us consider the vertices with $\deg_{G^{\prime}} \leq 3$
(note that such vertices must exist, since $E(G^{\prime}) < 2|G^{\prime}|$),
and let us copy the proof of Observation~\ref{observation}.
Observe firstly that there can't be any vertices with degree less than $2$,
since all faces in $G^{\prime}$ are cycles.
Secondly, note that if a vertex $v$ of degree $2$ exists,
then it must lie on the boundary of a $C_{3}$ face $v,v_{1},v_{2}$
and a $C_{6}$ face $v,v_{1},v_{3},v_{4},v_{5},v_{2}$,
as shown in Figure~\ref{trihexdeg2},
in which case $v_{1},v_{3},v_{4},v_{5},v_{2}$ would be a cycle of size $5$,
which is a contradiction.
Finally, note that if $\deg(v) = 3$ and $\Gamma(v) = \{v_{1},v_{2},v_{3}\}$,
then it is impossible for all the edges $vv_{1}$, $vv_{2}$ and $vv_{3}$ to be in a $C_{3}$ face
without having two adjacent such faces,
which would also be a contradiction. \\
\\
\textbf{Part IV} 
In this final part of the proof,
we shall deal with the remaining case when
$e(G^{\prime}) \leq 2|G^{\prime}| - 5$.
The proof will follow from a very careful examination of the impact of vertices of degree $2$ and $3$ in $G^{\prime}$.
We shall look first at vertices of degree $3$,
obtaining equation~(\ref{d3eqn}).
We shall then consider vertices of degree $2$,
producing (via equations~(\ref{case ii words}) and~(\ref{case ii symbols}) and Claim~\ref{f31 claim})
Claim~\ref{2k-5 claim} and equation~(\ref{importanteqn}).
Finally,
by plugging these newly derived inequalities into our `standard' framework
(i.e.~the fact that $2e(G)$ equals the sum of all face sizes,
plus Euler's formula),
we obtain the result.

Let $d_{i}^{\prime}$ denote the number of vertices of degree $i$ in $G^{\prime}$,
and recall that there can't be any vertices with degree less than $2$,
since all faces in $G^{\prime}$ are cycles.
Note that we must then have 
\begin{eqnarray}
2d_{2}^{\prime} + d_{3}^{\prime} & \geq & 
\left\{ \begin{array}{ll}
10 & \textrm{if $e(G^{\prime}) = 2|G^{\prime}| - 5$} \\
12 & \textrm{if $e(G^{\prime}) \leq 2|G^{\prime}| - 6$}
\end{array} \right. \label{d2d3}
\end{eqnarray}
since $2d_{2}^{\prime} + d_{3}^{\prime} \geq \sum_{i \geq 2} (4-i) d_{i}^{\prime} = 4|G^{\prime}| - 2e(G^{\prime})$.

Consider a vertex $v$ of degree $3$ in $G^{\prime}$,
and, in particular, consider the three faces that it lies on.
Note that if any of these is a $C_{3}$ or is a $C_{4}$ that was formed from two adjacent $C_{3}$'s in $G$,
then neither of the others can have size less than $6$
(this follows from Claims~\ref{triquadrilateralclaim} and~\ref{quadrilateralquadrilateral}
and the recollection that $G^{\prime}$ contains no adjacent $C_{3}$ faces).
Hence, it follows that either 
(i) none of the three edges incident to $v$ were in any $C_{3}$ faces in $G$ 
or (ii) there is at least one edge incident to $v$ that was in two faces of size at least $6$ in $G$.

Thus, we have
\begin{eqnarray}
2 \left( \frac{1}{3} (\# \textrm{ edges not in a $C_{3}$ face in $G$}) \phantom{wwwwwww} \right. & & \nonumber \\
+ \# \textrm{ edges in two faces of size at least $6$ in $G$}\Big) & \geq & d_{3}^{\prime} \nonumber
\end{eqnarray}
(where the multiplication of the left-hand side by $2$ is to allow for possible double-counting of an edge).
If we let $e_{i}$ denote the number of edges of $G$ that are in at least one $C_{i}$ face,
then this corresponds to
\begin{equation}
2 \left( \frac{1}{3} (e(G) - e_{3}) + (e(G) - e_{3} - e_{4}) \right) \geq d_{3}^{\prime}. \label{d3eqn}
\end{equation}

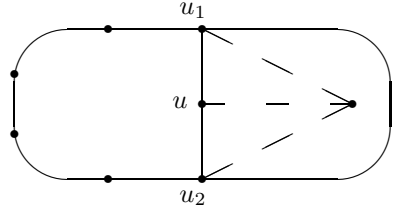
\begin{figure} [ht]
\setlength{\unitlength}{1cm}
\begin{picture}(20,2.8) (-6,-1.3)
\put(2,0){\line(-2,1){0.4}}
\put(1.2,0.4){\line(-2,1){0.4}}
\put(0.4,0.8){\line(-2,1){0.4}}
\put(2,0){\line(-2,-1){0.4}}
\put(1.2,-0.4){\line(-2,-1){0.4}}
\put(0.4,-0.8){\line(-2,-1){0.4}}
\put(2,0){\line(-1,0){0.3}}
\put(1.15,0){\line(-1,0){0.3}}
\put(0.3,0){\line(-1,0){0.3}}

\put(0,-1){\line(0,1){2}}
\put(0,0){\oval(5,2)[r]}
\put(0,-1){\circle*{0.1}}
\put(0,1){\circle*{0.1}}
\put(0,0){\circle*{0.1}}
\put(2,0){\circle*{0.1}}
\put(-0.3,1.2){$u_{1}$}
\put(-0.3,-1.3){$u_{2}$}
\put(-0.4,-0.05){$u$}
\put(0,0){\oval(5,2)[l]}

\put(-1.25,1){\circle*{0.1}}
\put(-1.25,-1){\circle*{0.1}}
\put(-2.5,0.4){\circle*{0.1}}
\put(-2.5,-0.4){\circle*{0.1}}
\end{picture}

\caption{The vertex $u$ in case (i).} \label{d2tri}
\end{figure}

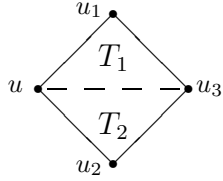
\begin{figure} [ht]
\setlength{\unitlength}{1cm}
\begin{picture}(20,2.4) (-6,-1.1)
\put(1,0){\line(-1,1){1}}
\put(1,0){\line(-1,-1){1}}
\put(-1,0){\line(1,1){1}}
\put(-1,0){\line(1,-1){1}}

\put(0,-1){\circle*{0.1}}
\put(0,1){\circle*{0.1}}
\put(-1,0){\circle*{0.1}}
\put(1,0){\circle*{0.1}}

\put(-1.4,-0.05){$u$}
\put(-0.5,1){$u_{1}$}
\put(-0.5,-1.1){$u_{2}$}
\put(1.1,-0.05){$u_{3}$}

\put(-0.875,0){\line(1,0){0.25}}
\put(-0.375,0){\line(1,0){0.25}}
\put(0.125,0){\line(1,0){0.25}}
\put(0.625,0){\line(1,0){0.25}}

\put(-0.2,0.3){\large{$T_{1}$}}
\put(-0.2,-0.6){\large{$T_{2}$}}

\end{picture}

\caption{The vertex $u$ in case (ii).} \label{d2quadrilateral}
\end{figure}

Now consider a vertex $u$ of degree $2$ in $G^{\prime}$,
and let $\Gamma_{G^{\prime}}(u) = \{u_{1},u_{2}\}$.
It must be that either
(i) $u$ lies on a $C_{3}$ face $uu_{1}u_{2}$ (formed from the outside of a $K_{4}$ in $G$),
as shown in Figure~\ref{d2tri},
or (ii) $u$ lies on a $C_{4}$ face $uu_{1}u_{3}u_{2}$ formed from two adjacent $C_{3}$ faces $T_{1} = uu_{1}u_{3}$ and $T_{2} = uu_{2}u_{3}$,
as shown in Figure~\ref{d2quadrilateral}.

In case (i), it will be important to note that the other face containing $u$ must have size at least $7$.
In case (ii), it will be important to note that the other faces bordering $T_{1}$ and $T_{2}$ in $G$ cannot have been $C_{3}$'s themselves
(since only the edge $uu_{3}$ was deleted when $G^{\prime}$ was constructed),
and so $T_{1}$ and $T_{2}$ are two $C_{3}$ faces in $G$
that both have exactly one edge ($uu_{3}$) that is also in another $C_{3}$ face in $G$.

Before proceeding,
let us also take care to allow for any possible multiple-counting.
For case (i),
note that it is possible for a face of size $i \geq 7$ to be adjacent in $G^{\prime}$
to $i-6$ faces of size $3$;
for case (ii),
note that both vertices $u$ and $u_{3}$ in Figure~\ref{d2quadrilateral} can have degree $2$ in $G^{\prime}$.

Thus, allowing for this multiple-counting,
we obtain
\begin{eqnarray}
2 \left( \frac{1}{2} (\# \textrm{ $C_{3}$ faces in $G$ that have exactly one} \phantom{www} \right. & & \nonumber \\
\textrm{ edge that is also in another $C_{3}$ face in $G$}) \Big) & & \nonumber \\ 
+ \# \textrm{ $C_{7}$ faces in $G$} + 2(\# \textrm{ $C_{8}$ faces in $G$}) & & \nonumber \\
+ 3(\# \textrm{ $C_{9}$ faces in $G$}) + \ldots & \geq & d_{2}^{\prime}. \label{case ii words}
\end{eqnarray}

Let $f_{i}$ denote the number of $C_{i}$ faces in $G$,
and let $f_{3,1}$ denote the number of $C_{3}$ faces in $G$ that have exactly one edge that is also in another $C_{3}$ face in $G$.
Then (\ref{case ii words}) translates to
\begin{equation}
f_{3,1} + \sum_{i \geq 7} (i-6)f_{i} \geq d_{2}^{\prime}. \label{case ii symbols}
\end{equation}

\begin{Claim} \label{2k-5 claim}
If $e(G^{\prime}) = 2|G^{\prime}|-5$, then $\sum_{i \geq 7} (i-6)f_{i} \geq d_{2}^{\prime} - 2$.
\end{Claim}
\textbf{Proof}
Note that, since $G \not\supset C_{5}$, $f_{3,1}$ must be equal to $2l$,
where $l$ is precisely the number of $C_{4}$ faces of $G^{\prime}$ that were formed from the remains of two adjacent $C_{3}$ faces of $G$.

Let $e_{4,1}^{\prime}$ denote the number of edges of $G^{\prime}$
that are in exactly one $C_{4}$ face in $G^{\prime}$.
Then,
by Claim~\ref{quadrilateralquadrilateral},
we have $e_{4,1}^{\prime} \geq 4l$.
Note that $e_{4}^{\prime} - e_{4,1}^{\prime}$
will be precisely the number of edges of $G^{\prime}$
that are in exactly two $C_{4}$ faces in $G^{\prime}$,
and thus
\begin{eqnarray}
f_{4}^{\prime} & = & \frac{ e_{4,1}^{\prime} + 2 ( e_{4}^{\prime} - e_{4,1}^{\prime} ) }{4} \nonumber \\
& = & \frac{ 2  e_{4}^{\prime} - e_{4,1}^{\prime} }{4} \nonumber \\
& \leq & \frac{2e_{4}^{\prime} - 4l}{4} \nonumber \\
& = & \frac{e_{4}^{\prime}}{2} - l. \label{f4}
\end{eqnarray}

We know that $2e(G^{\prime}) \geq 3f_{3}^{\prime} + 4f_{4}^{\prime} + 6(f^{\prime} - f_{3}^{\prime} - f_{4}^{\prime})$,
and so
\begin{eqnarray*}
6f^{\prime} & \leq & 2e(G^{\prime}) + 3f_{3}^{\prime} + 2f_{4}^{\prime} \\
& \leq & 2e(G^{\prime}) + e_{3}^{\prime} + e_{4}^{\prime} - 2l \textrm{ using (\ref{f3prime}) and (\ref{f4})} \\
& \leq & 2e(G^{\prime}) + e_{3}^{\prime} + e(G^{\prime}) - e_{3}^{\prime} - 2l \textrm{ by Claim~\ref{triquadrilateralclaim}} \\
& = & 3e(G^{\prime}) - 2l.
\end{eqnarray*}
Hence, $f^{\prime} \leq \frac{1}{2} e(G^{\prime}) - \frac{l}{3}$.

Euler's formula $e(G^{\prime}) = |G^{\prime}| - 2 + f^{\prime}$ then gives $e(G^{\prime}) \leq 2 \left( |G^{\prime}| - \left(2 + \frac{l}{3} \right) \right)$.
If $e(G^{\prime}) = 2|G^{\prime}| - 5$,
this implies that $l \leq 1$ (since $l$ is an integer),
and hence that $f_{3,1} \leq 2$.
The result then follows from (\ref{case ii symbols}).
\textbf{End of Proof of Claim~\ref{2k-5 claim}} 

\begin{Claim} \label{f31 claim}
$f_{3,1} \leq 2e_{3} - 4f_{3}$.
\end{Claim}
\textbf{Proof}
Note that it is impossible for all three edges in a $C_{3}$ face in $G$ to also be in another $C_{3}$ face in $G$ without creating a $C_{5}$
unless $G=C_{3}$ or $G=K_{4}$,
both of which may be ignored since we know that $k>4$.
Thus,
\begin{eqnarray*}
f_{3} & = & \frac{1}{3} ( \# \textrm{ edges in at least one $C_{3}$ face in $G$ } + \# \textrm{ edges in two $C_{3}$ faces in $G$}) \\
& \leq & \frac{1}{3} \left( e_{3} + \frac{1}{2} ( f_{3,1} + 2(f_{3} - f_{3,1})) \right) \\
& = & \frac{1}{3}e_{3} + \frac{1}{3}f_{3} - \frac{1}{6}f_{3,1},
\end{eqnarray*}
from which the claim follows.
\textbf{End of Proof of Claim~\ref{f31 claim}} \\

By (\ref{case ii symbols}) and Claim~\ref{f31 claim}, we have
\begin{equation} \label{importanteqn}
2e_{3} - 4f_{3} + \sum_{i \geq 7} (i-6)f_{i} \geq d_{2}^{\prime}.
\end{equation}

We shall now use (\ref{d3eqn}), (\ref{importanteqn}) and Claim~\ref{2k-5 claim} to obtain our result.
Letting $f$ denote $\sum_{i} f_{i}$,
we know that $2e(G) = 3f_{3} + 4f_{4} + 6(f-f_{3}-f_{4}) + \sum_{i \geq 7}(i-6)f_{i}$,
so we have
\begin{eqnarray*}
6f & = & 2e(G) + 3f_{3} + 2f_{4} - \sum_{i \geq 7}(i-6)f_{i} \\
& \leq & 2e(G) + 3f_{3} + e_{4}  - \sum_{i \geq 7}(i-6)f_{i} \\
& = & 2e(G) + \frac{3}{2}e_{3} + e_{4}  - \sum_{i \geq 7}(i-6)f_{i} - \frac{3}{4}(2e_{3} - 4f_{3}) \\
& = & 2e(G) + \frac{3}{2}e_{3} + e(G) - e_{3}  - \sum_{i \geq 7}(i-6)f_{i} - \frac{3}{4}(2e_{3} - 4f_{3}) - (e(G) - e_{3} -e_{4}) \\
& = & 3e(G) + \frac{1}{2}e_{3} - \sum_{i \geq 7}(i-6)f_{i} - \frac{3}{4}(2e_{3} - 4f_{3}) - (e(G) - e_{3} -e_{4}) \\
& = & \frac{7}{2}e(G) - \sum_{i \geq 7}(i-6)f_{i} - \frac{3}{4}(2e_{3} - 4f_{3}) - (e(G) - e_{3} -e_{4}) - \frac{1}{2}(e(G)-e_{3}) \\
& = & \frac{7}{2}e(G) - \frac{3}{4} \left( 2e_{3} - 4f_{3} + \sum_{i \geq 7}(i-6)f_{i} \right) 
- \left( \frac{1}{3} (e(G) - e_{3}) + (e(G) - e_{3} -e_{4}) \right) \\
& & - \frac{1}{4} \sum_{i \geq 7}(i-6)f_{i} - \frac{1}{6}(e(G)-e_{3}) \\
& \leq & \frac{7}{2} e(G) - \frac{3}{4} d_{2}^{\prime} - \frac{1}{2} d_{3}^{\prime} - \frac{1}{4} \sum_{i \geq 7}(i-6)f_{i} \textrm{ by (\ref{d3eqn}) and (\ref{importanteqn}).}
\end{eqnarray*}

Thus, $f \leq \frac{7}{12} e(G) - \frac{1}{8} d_{2}^{\prime} - \frac{1}{12} d_{3}^{\prime} - \frac{1}{24} \sum_{i \geq 7}(i-6)f_{i}$,
and so Euler's formula $e(G) = k-2+f$ then implies that
\begin{eqnarray*}
e(G) & \leq & \frac{12}{5}(k - (2 + \frac{1}{8}d_{2}^{\prime} + \frac{1}{12}d_{3}^{\prime} + \frac{1}{24} \sum_{i \geq 7}(i-6)f_{i})) \\
& = & \frac{12k - (24 + \frac{3}{2} d_{2}^{\prime} + d_{3}^{\prime} + \frac{1}{2} \sum_{i \geq 7}(i-6)f_{i})}{5}.
\end{eqnarray*}

It now only remains to show that $\frac{3}{2} d_{2}^{\prime} + d_{3}^{\prime} + \frac{1}{2} \sum_{i \geq 7}(i-6)f_{i} \geq 9$.
If $e(G^{\prime}) \leq 2|G^{\prime}| - 6$, then $2d_{2}^{\prime} + d_{3}^{\prime} \geq 12$
(from (\ref{d2d3})),
and so certainly $\frac{3}{2} d_{2}^{\prime} + d_{3}^{\prime} \geq 9$.
If $e(G^{\prime}) = 2 |G^{\prime}| - 5$,
then recall (from Claim~\ref{2k-5 claim}) that we have
$\sum_{i \geq 7} (i-6)f_{i} \geq d_{2}^{\prime} - 2$,
and hence $\frac{3}{2} d_{2}^{\prime} + d_{3}^{\prime} + \frac{1}{2} \sum_{i \geq 7}(i-6)f_{i} \geq 2 d_{2}^{\prime} + d_{3}^{\prime} - 1 \geq 9$
(using (\ref{d2d3})).
\phantom{qwerty}
\setlength{\unitlength}{0.25cm}
\begin{picture}(1,1)
\put(0,0){\line(1,0){1}}
\put(0,0){\line(0,1){1}}
\put(1,1){\line(-1,0){1}}
\put(1,1){\line(0,-1){1}}
\end{picture} \\

\section{$\mathbf{C_{5}}$ --- extremal graphs} \label{examplessection}

In the previous section,
we showed that $\ex(n,C_{5}) \leq \frac{12n-33}{5}$ for all $n \geq 11$.
In this section
we shall now complete matters (in Theorem~\ref{endthm}) by demonstrating that this bound is tight.

\begin{figure} [ht]
\setlength{\unitlength}{1cm}
\begin{picture}(20,2.4) (-5.5,-1.1)
\put(1,0){\line(-1,1){1}}
\put(1,0){\line(-1,-1){1}}
\put(-1,0){\line(1,1){1}}
\put(-1,0){\line(1,-1){1}}
\put(0,0){\oval(4,2)[r]}
\put(0,-1){\circle*{0.1}}
\put(0,1){\circle*{0.1}}
\put(-1,0){\circle*{0.1}}
\put(1,0){\circle*{0.1}}
\put(-1.4,-0.05){$v_{1}$}
\put(-0.45,1){$v_{2}$}
\put(-0.4,-1.1){$v_{5}$}
\put(1.1,-0.05){$v_{6}$}

\put(0,-1){\line(0,1){2}}
\put(-1,0){\line(3,1){1}}
\put(-1,0){\line(3,-1){1}}
\put(1,0){\line(-3,1){1}}
\put(1,0){\line(-3,-1){1}}
\put(0,0.333333){\circle*{0.1}}
\put(0,-0.333333){\circle*{0.1}}
\put(-0.35,0.45){$v_{3}$}
\put(-0.35,-0.15){$v_{4}$}

\end{picture}
\caption{The triangulation $T_{6}$.} \label{T6tri}
\end{figure}
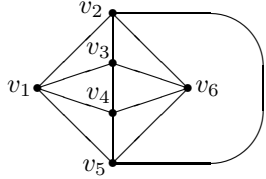

We start with a lemma that will prove useful:

\begin{Lemma} \label{endlemma}
For infinitely many values of $k$,
there exists a plane triangulation $T_{k}$ with vertex set $\{v_{1}, v_{2}, \ldots, v_{k} \}$ satisfying
(i)~$\deg (v_{i}) = 4$ for $i \leq 6$,
(ii)~$\deg (v_{i}) = 6$ for $i > 6$,
and (iii)~$E(T_{k}) \supset \{ v_{1}v_{2}, v_{3}v_{4}, v_{5}v_{6} \}$.
\end{Lemma}
\textbf{Proof}
Note first that the triangulation $T_{6}$ shown in Figure~\ref{T6tri} certainly satisfies the conditions for $k=6$.

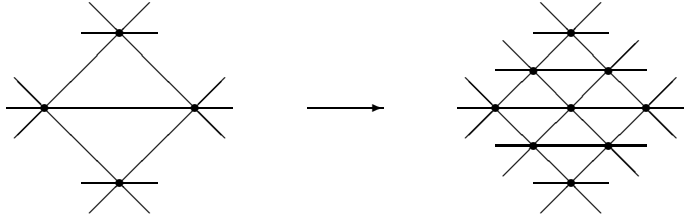
\begin{figure} [ht]
\setlength{\unitlength}{1cm}
\begin{picture}(20,3) (-3,-1.5)

\put(1,0){\line(-1,1){1.4}}
\put(1,0){\line(-1,-1){1.4}}
\put(-1,0){\line(1,1){1.4}}
\put(-1,0){\line(1,-1){1.4}}
\put(0,-1){\circle*{0.1}}
\put(0,1){\circle*{0.1}}
\put(-1,0){\circle*{0.1}}
\put(1,0){\circle*{0.1}}
\put(-1.5,0){\line(1,0){3}}
\put(-1,0){\line(-1,1){0.4}}
\put(-1,0){\line(-1,-1){0.4}}
\put(1,0){\line(1,1){0.4}}
\put(1,0){\line(1,-1){0.4}}

\put(-0.5,1){\line(1,0){1}}
\put(-0.5,-1){\line(1,0){1}}

\put(2.5,0){\vector(1,0){1}}

\put(6,0){

\put(1,0){\line(-1,1){1.4}}
\put(1,0){\line(-1,-1){1.4}}
\put(-1,0){\line(1,1){1.4}}
\put(-1,0){\line(1,-1){1.4}}
\put(0,-1){\circle*{0.1}}
\put(0,1){\circle*{0.1}}
\put(-1,0){\circle*{0.1}}
\put(1,0){\circle*{0.1}}
\put(-1.5,0){\line(1,0){3}}
\put(-1,0){\line(-1,1){0.4}}
\put(-1,0){\line(-1,-1){0.4}}
\put(1,0){\line(1,1){0.4}}
\put(1,0){\line(1,-1){0.4}}

\put(-0.5,1){\line(1,0){1}}
\put(-0.5,-1){\line(1,0){1}}

}

\put(5.5,0.5){\circle*{0.1}}
\put(5.5,-0.5){\circle*{0.1}}
\put(6.5,0.5){\circle*{0.1}}
\put(6.5,-0.5){\circle*{0.1}}
\put(6,0){\circle*{0.1}}

\put(5,0.5){\line(1,0){2}}
\put(5,-0.5){\line(1,0){2}}
\put(5.1,0.9){\line(1,-1){1.8}}
\put(5.1,-0.9){\line(1,1){1.8}}

\end{picture}
\caption{Constructing a larger triangulation.} \label{largertri}
\end{figure}

We now proceed inductively.
Given a triangulation satisfying the conditions,
let us construct a larger triangulation by subdividing all the edges 
and inserting triangles between the new vertices,
as shown in Figure~\ref{largertri}.

\begin{figure} [ht]
\setlength{\unitlength}{1cm}
\begin{picture}(20,3) (-3,-1.5)
\put(1,0){\line(-1,1){1.4}}
\put(1,0){\line(-1,-1){1.4}}
\put(-1,0){\line(1,1){1.4}}
\put(-1,0){\line(1,-1){1.4}}
\put(0,-1){\circle*{0.1}}
\put(0,1){\circle*{0.1}}
\put(-1,0){\circle*{0.1}}
\put(1,0){\circle*{0.1}}
\put(-1.5,0){\line(1,0){3}}

\put(-0.5,1){\line(1,0){1}}
\put(-0.5,-1){\line(1,0){1}}

\put(-0.5,0.5){\circle*{0.1}}
\put(-0.5,-0.5){\circle*{0.1}}
\put(0.5,0.5){\circle*{0.1}}
\put(0.5,-0.5){\circle*{0.1}}
\put(0,0){\circle*{0.1}}

\put(-1,0.5){\line(1,0){2}}
\put(-1,-0.5){\line(1,0){2}}
\put(-0.9,0.9){\line(1,-1){1.8}}
\put(-0.9,-0.9){\line(1,1){1.8}}

\put(2.5,0){\vector(1,0){1}}

\put(6,0){

\put(1,0){\line(-1,1){1.4}}
\put(1,0){\line(-1,-1){1.4}}
\put(-1,0){\line(1,1){1.4}}
\put(-1,0){\line(1,-1){1.4}}
\put(0,-1){\circle*{0.1}}
\put(0,1){\circle*{0.1}}
\put(-1,0){\circle*{0.1}}
\put(1,0){\circle*{0.1}}
\put(-1.5,0){\line(1,0){3}}

\put(-0.5,1){\line(1,0){1}}
\put(-0.5,-1){\line(1,0){1}}

\put(-0.5,0.5){\circle*{0.1}}
\put(-0.5,-0.5){\circle*{0.1}}
\put(0.5,0.5){\circle*{0.1}}
\put(0.5,-0.5){\circle*{0.1}}

\put(-1,0.5){\line(1,0){2}}
\put(-1,-0.5){\line(1,0){2}}
\put(-0.9,0.9){\line(1,-1){0.4}}
\put(-0.9,-0.9){\line(1,1){0.4}}

}

\put(6.5,0.5){\line(1,1){0.4}}
\put(6.5,-0.5){\line(1,-1){0.4}}

\put(6,0.4){\line(5,1){0.5}}
\put(6,0.4){\line(-5,1){0.5}}
\put(6,0.4){\line(5,-2){1}}
\put(6,0.4){\line(1,-4){0.1}}
\put(6,0.4){\line(-1,-4){0.1}}
\put(6,0.4){\line(-5,-2){1}}
\put(6,-0.4){\line(5,-1){0.5}}
\put(6,-0.4){\line(-5,-1){0.5}}
\put(6,-0.4){\line(5,2){1}}
\put(6,-0.4){\line(1,4){0.1}}
\put(6,-0.4){\line(-1,4){0.1}}
\put(6,-0.4){\line(-5,2){1}}

\put(5.9,0){\circle*{0.1}}
\put(6,-0.4){\circle*{0.1}}
\put(6,0.4){\circle*{0.1}}
\put(6.1,0){\circle*{0.1}}

\put(-1.3,-0.3){$v_{1}$}
\put(0.9,-0.3){$v_{2}$}

\put(4.6,-0.35){$v_{1}^{\prime}$}
\put(6.9,-0.35){$v_{2}^{\prime}$}

\put(5.55,0.05){\small{$v_{1}$}}
\put(6.15,0.05){\small{$v_{2}$}}

\end{picture}
\caption{Modifying the triangulation to satisfy condition (iii).} \label{modifytri}
\end{figure}
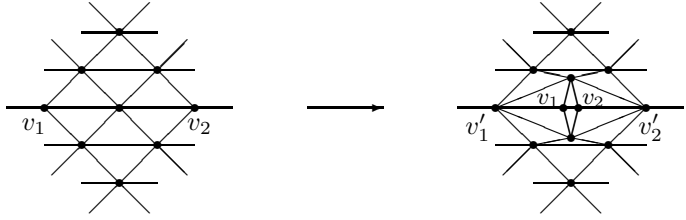

It can be observed that conditions (i) and (ii) will also be satisfied by the new triangulation,
but (due to the subdividing of edges) not condition (iii).
However,
we may then simply modify the new triangulation into one that does satisfy all three conditions 
by applying the local transformation shown in Figure~\ref{modifytri} 
(which includes some relabelling of the vertices)
at the relevant three places.
\phantom{qwerty}
\setlength{\unitlength}{0.25cm}
\begin{picture}(1,1)
\put(0,0){\line(1,0){1}}
\put(0,0){\line(0,1){1}}
\put(1,1){\line(-1,0){1}}
\put(1,1){\line(0,-1){1}}
\end{picture} \\

We may now proceed with our construction of extremal graphs:

\begin{Theorem} \label{endthm}
There exist infinitely many values of $n$ for which
$\ex(n,C_{5}) = \frac{12n-33}{5}$.
\end{Theorem}
\textbf{Proof} 
We shall use the triangulations guaranteed by Lemma~\ref{endlemma}
to construct specific plane graphs with $15k+9$ vertices and $\frac{12(15k+9)-33}{5}$ edges,
and we shall then show that these graphs are $C_{5}$-free.

\begin{figure} [ht]
\setlength{\unitlength}{1cm}
\begin{picture}(20,2.3) (-8,-1.15)

\put(-5,-1){\circle*{0.2}}
\put(-5,1){\circle*{0.2}}
\put(-5,-1){\line(0,1){2}}

\put(-3.5,0){\vector(1,0){1}}

\put(1,0){\line(-1,1){1}}
\put(1,0){\line(-1,-1){1}}
\put(-1,0){\line(1,1){1}}
\put(-1,0){\line(1,-1){1}}

\put(0,-1){\circle*{0.2}}
\put(0,1){\circle*{0.2}}
\put(-1,0){\circle*{0.1}}
\put(1,0){\circle*{0.1}}

\put(0,0.4){\circle*{0.1}}
\put(0,-0.4){\circle*{0.1}}

\put(0,-1){\line(0,1){0.6}}
\put(0,1){\line(0,-1){0.6}}

\put(0,0.4){\line(1,-1){0.4}}
\put(0,0.4){\line(-1,-1){0.4}}
\put(0,-0.4){\line(1,1){0.4}}
\put(0,-0.4){\line(-1,1){0.4}}
\put(-0.4,0){\line(1,0){0.8}}
\put(-0.4,0){\circle*{0.1}}
\put(0.4,0){\circle*{0.1}}

\end{picture}
\caption{Replacing an edge with a `diamond-holder'.} \label{diamondholder}
\end{figure}
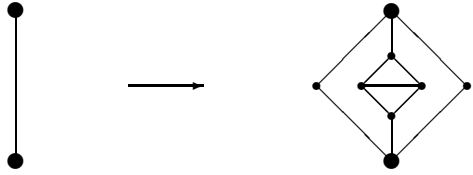

\begin{figure} [ht]
\setlength{\unitlength}{1cm}
\begin{picture}(20,6.85) (-9,1.4)
\put(-6.5,0){

\put(0,2){\line(0,1){5.75}}

\put(1.5,3){\line(-3,2){1.5}}
\put(1.5,5){\line(-3,2){1.5}}
\put(-1.5,3){\line(3,2){1.5}}
\put(-1.5,5){\line(3,2){1.5}}

\put(1.5,3){\line(0,1){4}}
\put(-1.5,3){\line(0,1){4}}

\put(1.5,5){\oval(1,4)[r]}
\put(-1.5,5){\oval(1,4)[l]}

\put(2,4.75){\oval(1,5.5)[r]}
\put(-2,4.75){\oval(1,5.5)[l]}

\put(2,7){\oval(1,1)[tl]}
\put(-2,7){\oval(1,1)[tr]}

\put(-2,2){\line(1,0){4}}

\put(0.5,7.75){\oval(1,1)[tl]}
\put(0.5,2){\oval(1,1)[bl]}
\put(0.5,4.875){\oval(5,6.75)[r]}

\multiput(0,0)(0,2){3}{

\put(0,2.5){\circle*{0.2}}
\put(0,2){\circle*{0.2}}
\put(0,3.5){\circle*{0.2}}
\put(1.5,3){\circle*{0.2}}
\put(-1.5,3){\circle*{0.2}}

\put(1.5,3){\line(-3,1){1.5}}
\put(1.5,3){\line(-3,-1){1.5}}
\put(1.5,3){\line(-3,-2){1.5}}
\put(-1.5,3){\line(3,1){1.5}}
\put(-1.5,3){\line(3,-1){1.5}}
\put(-1.5,3){\line(3,-2){1.5}}

}

}

\put(-8.5,1.4){\Large{$T_{15}$}}

\put(-3.25,5){\vector(1,0){0.5}}

\put(1.5,3){\line(-3,2){1.5}}
\put(1.5,5){\line(-3,2){1.5}}
\put(-1.5,3){\line(3,2){1.5}}
\put(-1.5,5){\line(3,2){1.5}}

\put(1.5,3){\line(0,1){4}}
\put(-1.5,3){\line(0,1){4}}

\put(1.5,5){\oval(1,4)[r]}
\put(-1.5,5){\oval(1,4)[l]}

\put(2,4.75){\oval(1,5.5)[r]}
\put(-2,4.75){\oval(1,5.5)[l]}

\put(2,7){\oval(1,1)[tl]}
\put(-2,7){\oval(1,1)[tr]}

\put(-2,2){\line(1,0){4}}

\put(0.5,7.75){\oval(1,1)[tl]}
\put(0.5,2){\oval(1,1)[bl]}
\put(0.5,4.875){\oval(5,6.75)[r]}

\multiput(0,0)(0,2){3}{

\put(0,2.5){\circle*{0.2}}
\put(0,2){\circle*{0.2}}
\put(0,3.5){\circle*{0.2}}
\put(1.5,3){\circle*{0.2}}
\put(-1.5,3){\circle*{0.2}}

\put(1.5,3){\line(-3,1){1.5}}
\put(1.5,3){\line(-3,-1){1.5}}
\put(1.5,3){\line(-3,-2){1.5}}
\put(-1.5,3){\line(3,1){1.5}}
\put(-1.5,3){\line(3,-1){1.5}}
\put(-1.5,3){\line(3,-2){1.5}}

}

\multiput(0,0)(0,2){3}{

\put(-0.5,3){\line(1,0){1}}

\put(0,2.75){\line(2,1){0.5}}
\put(0,3.25){\line(2,-1){0.5}}
\put(0.5,3){\circle*{0.1}}
\put(0,3.5){\line(2,-1){1}}
\put(0,2.5){\line(2,1){1}}
\put(1,3){\circle*{0.1}}

\put(0,2.75){\line(-2,1){0.5}}
\put(0,3.25){\line(-2,-1){0.5}}
\put(-0.5,3){\circle*{0.1}}
\put(0,3.5){\line(-2,-1){1}}
\put(0,2.5){\line(-2,1){1}}
\put(-1,3){\circle*{0.1}}

\put(0,3.25){\circle*{0.1}}
\put(0,2.75){\circle*{0.1}}
\put(0,2.5){\circle*{0.2}}

\put(0,2){\circle*{0.2}}

\put(0,2.25){\circle*{0.1}}

\put(0,3.5){\circle*{0.2}}
\put(0,3.75){\circle*{0.1}}

\put(0,2){\line(0,1){0.75}}

\put(0.75,2.5){\circle*{0.1}}
\put(0.75,2.75){\circle*{0.1}}
\put(0.75,3.25){\circle*{0.1}}

\put(-0.75,2.5){\circle*{0.1}}
\put(-0.75,2.75){\circle*{0.1}}
\put(-0.75,3.25){\circle*{0.1}}

}

\put(0,3.25){\line(0,1){1.25}}
\put(0,5.25){\line(0,1){1.25}}
\put(0,7.25){\line(0,1){0.5}}

\put(1.5,4){\circle*{0.1}}
\put(1.5,6){\circle*{0.1}}
\put(2,5){\circle*{0.1}}
\put(2.5,5){\circle*{0.1}}
\put(-1.5,4){\circle*{0.1}}
\put(-1.5,6){\circle*{0.1}}
\put(-2,5){\circle*{0.1}}
\put(-2.5,5){\circle*{0.1}}

\put(0.75,3.5){\circle*{0.1}}
\put(0.75,5.5){\circle*{0.1}}

\put(-0.75,3.5){\circle*{0.1}}
\put(-0.75,5.5){\circle*{0.1}}

\put(-1.8,1.4){\Large{$G$}}

\end{picture}
\caption{Constructing a new graph $G$ from a triangulation $T_{15}$.} \label{GfromT15}
\end{figure}
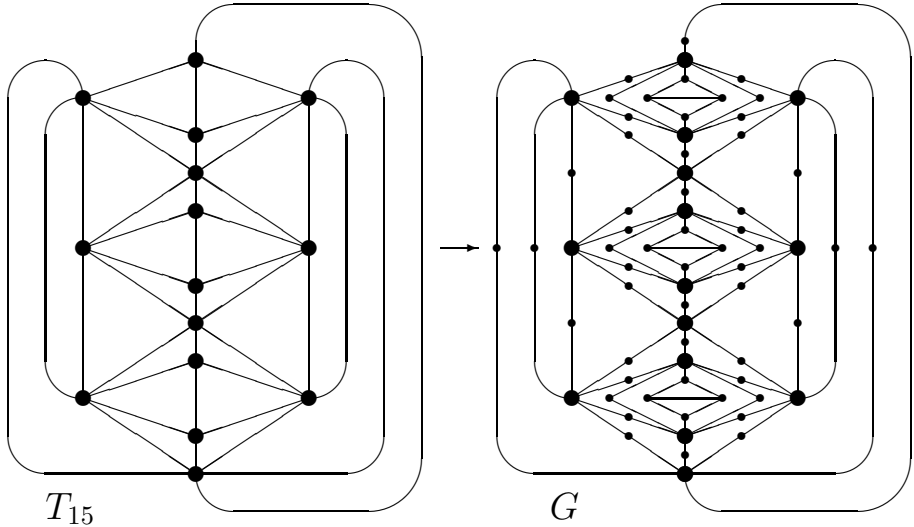

Let $T_{k}$ be a plane triangulation satisfying the conditions of Lemma~\ref{endlemma},
and let $E^{*}$ denote the set of edges $\{v_{1}v_{2}, v_{3}v_{4}, v_{5}v_{6}\}$.
Let us now construct a new plane graph $G$ by 
(i)~subdividing all edges in $E(T_{k}) \setminus E^{*}$
and (ii)~replacing all edges in $E^{*}$ with the `diamond-holder' structure shown in Figure~\ref{diamondholder}.
An example of this is shown in Figure~\ref{GfromT15}.

Note that the $k$ original vertices 
(shown as large vertices in Figure~\ref{GfromT15})
will now all have degree $6$
(including the six that only had degree $4$ in $T_{k}$),
and so each can be thought of as forming the centre of a `star' of size six (i.e.~$K_{1,6}$).
Hence,
$G$ consists of $k$ of these stars,
together with the three `diamonds' 
(i.e.~$K_{4}$ minus an edge)
from the centres of the diamond-holders.

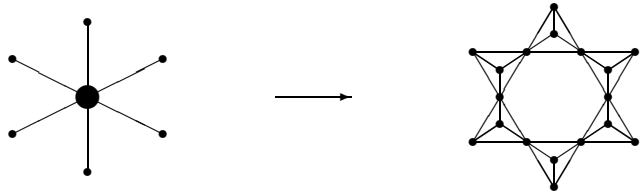
\begin{figure} [ht]
\setlength{\unitlength}{1cm}
\begin{picture}(20,2.4) (-3,-0.2)

\put(0,0){\line(0,1){2}}
\put(-1,0.5){\line(2,1){2}}
\put(-1,1.5){\line(2,-1){2}}
\put(0,1){\circle*{0.3}}
\put(-1,0.5){\circle*{0.1}}
\put(-1,1.5){\circle*{0.1}}
\put(1,0.5){\circle*{0.1}}
\put(1,1.5){\circle*{0.1}}
\put(0,0){\circle*{0.1}}
\put(0,2){\circle*{0.1}}

\put(2.5,1){\vector(1,0){1}}

\put(5,0.2){

\put(0.12,0.2){\line(1,0){2.16}}
\put(0.12,0.2){\line(3,5){1.08}}
\put(2.28,0.2){\line(-3,5){1.08}}
\put(0.12,0.2){\circle*{0.1}}
\put(2.28,0.2){\circle*{0.1}}
\put(1.2,2){\circle*{0.1}}

\put(0.12,1.4){\line(1,0){2.16}}
\put(2.28,1.4){\line(-3,-5){1.08}}
\put(0.12,1.4){\line(3,-5){1.08}}
\put(1.2,-0.4){\circle*{0.1}}
\put(2.28,1.4){\circle*{0.1}}
\put(0.12,1.4){\circle*{0.1}}

\put(0.84,0.2){\circle*{0.1}}
\put(1.56,0.2){\circle*{0.1}}
\put(0.48,0.8){\circle*{0.1}}
\put(1.92,0.8){\circle*{0.1}}
\put(0.84,1.4){\circle*{0.1}}
\put(1.56,1.4){\circle*{0.1}}

\put(1.2,2){\line(0,-1){0.36}}
\put(1.2,1.64){\line(-3,-2){0.72}}
\put(1.2,1.64){\line(3,-2){0.72}}
\put(1.2,1.64){\circle*{0.1}}

\put(1.92,1.16){\line(0,-1){0.72}}
\put(0.48,1.16){\line(0,-1){0.72}}
\put(1.92,1.16){\circle*{0.1}}
\put(0.48,1.16){\circle*{0.1}}
\put(1.92,0.44){\circle*{0.1}}
\put(0.48,0.44){\circle*{0.1}}

\put(1.2,-0.4){\line(0,1){0.36}}
\put(1.2,-0.04){\line(-3,2){0.72}}
\put(1.2,-0.04){\line(3,2){0.72}}
\put(1.2,-0.04){\circle*{0.1}}

\put(1.92,1.16){\line(3,2){0.36}}
\put(0.48,1.16){\line(-3,2){0.36}}
\put(1.92,0.44){\line(3,-2){0.36}}
\put(0.48,0.44){\line(-3,-2){0.36}}

}

\end{picture}
\caption{Replacing a star with a `snowflake'.} \label{snowflake}
\end{figure}

Finally,
let us now replace each star with the `snowflake' structure shown in Figure~\ref{snowflake},
where the central vertex is replaced by a hexagon
and the edges are replaced by $K_{4}$'s.
Thus, the final graph $G^{*}$ will consist of $k$ snowflakes and three diamonds.

Observe that each diamond contains $4$ vertices,
$2$ of which will be shared with snowflakes,
and $5$ edges;
each snowflake contains $18$ vertices,
$6$ of which will be shared with other snowflakes/diamonds,
and $36$ edges.

Hence,
allowing for double-counting of the shared vertices,
we find that $G^{*}$ contains
$3(2 + \frac{2}{2}) + k (12 + \frac{6}{2}) = 15k+9$ vertcies
and $3 \times 5 + 36k = \frac{12(15k+9)-33}{5}$ edges,
thus achieving the bound of Theorem~\ref{mainthm}.

It now only remains to show that $G^{*}$ is $C_{5}$-free.

To see this,
let us first examine the earlier graph $G$.
Due to the subdividing of the edges of $T_{k}$
(and the resulting doubling of all cycle lengths),
note that $G$ can contain no cycles of length less than six
apart from those
(of length three and four)
contained entirely within a diamond-holder.

Now observe that any cycle in $G^{*}$ of length less than six must 
either induce a cycle of length less than six in $G$
or must be contained entirely within a single snowflake.

\begin{figure} [ht]
\setlength{\unitlength}{0.75cm}
\begin{picture}(20,9) (-8,-4.5)

\put(-0.5,3.5){\line(1,0){1}}
\put(-0.5,1.5){\line(1,0){1}}
\put(1,2.5){\line(-1,2){0.5}}
\put(1,2.5){\line(-1,-2){0.5}}
\put(-1,2.5){\line(1,2){0.5}}
\put(-1,2.5){\line(1,-2){0.5}}
\put(-0.5,3.5){\circle*{0.133333}}
\put(0.5,3.5){\circle*{0.133333}}
\put(-1,2.5){\circle*{0.133333}}
\put(1,2.5){\circle*{0.133333}}
\put(-0.5,1.5){\circle*{0.133333}}
\put(0.5,1.5){\circle*{0.133333}}

\put(0.5,1.5){\line(-1,-2){0.5}}
\put(0.5,1.5){\line(-1,-1){0.5}}
\put(0.5,1.5){\line(1,0){0.5}}
\put(1,2.5){\line(0,-1){1}}
\put(0.5,1.5){\line(2,-3){1}}
\put(1,1.5){\line(1,-3){0.5}}
\put(1,2.5){\line(1,-5){0.5}}
\put(0,0.5){\line(1,-1){0.5}}
\put(0.5,0){\circle*{0.133333}}
\put(1,1.5){\circle*{0.133333}}
\put(1.5,0){\circle*{0.133333}}

\put(0,1){\line(0,-1){0.5}}
\put(-0.5,0){\line(1,0){1}}
\put(0,1){\circle*{0.133333}}
\put(0,0.5){\circle*{0.133333}}

\put(-0.5,1.5){\line(1,-2){0.5}}
\put(-0.5,1.5){\line(1,-1){0.5}}
\put(-0.5,1.5){\line(-1,0){0.5}}
\put(-1,2.5){\line(0,-1){1}}
\put(-0.5,1.5){\line(-2,-3){1}}
\put(-1,1.5){\line(-1,-3){0.5}}
\put(-1,2.5){\line(-1,-5){0.5}}
\put(0,0.5){\line(-1,-1){0.5}}
\put(-0.5,0){\circle*{0.133333}}
\put(-1,1.5){\circle*{0.133333}}
\put(-1.5,0){\circle*{0.133333}}

\put(-0.5,-3.5){\line(1,0){1}}
\put(-0.5,-1.5){\line(1,0){1}}
\put(1,-2.5){\line(-1,-2){0.5}}
\put(1,-2.5){\line(-1,2){0.5}}
\put(-1,-2.5){\line(1,-2){0.5}}
\put(-1,-2.5){\line(1,2){0.5}}
\put(-0.5,-3.5){\circle*{0.133333}}
\put(0.5,-3.5){\circle*{0.133333}}
\put(-1,-2.5){\circle*{0.133333}}
\put(1,-2.5){\circle*{0.133333}}
\put(-0.5,-1.5){\circle*{0.133333}}
\put(0.5,-1.5){\circle*{0.133333}}

\put(0.5,-1.5){\line(-1,2){0.5}}
\put(0.5,-1.5){\line(-1,1){0.5}}
\put(0.5,-1.5){\line(1,0){0.5}}
\put(1,-2.5){\line(0,1){1}}
\put(0.5,-1.5){\line(2,3){1}}
\put(1,-1.5){\line(1,3){0.5}}
\put(1,-2.5){\line(1,5){0.5}}
\put(0,-0.5){\line(1,1){0.5}}
\put(1,-1.5){\circle*{0.133333}}

\put(0,-1){\line(0,1){0.5}}
\put(0,-1){\circle*{0.133333}}
\put(0,-0.5){\circle*{0.133333}}

\put(-0.5,-1.5){\line(1,2){0.5}}
\put(-0.5,-1.5){\line(1,1){0.5}}
\put(-0.5,-1.5){\line(-1,0){0.5}}
\put(-1,-2.5){\line(0,1){1}}
\put(-0.5,-1.5){\line(-2,3){1}}
\put(-1,-1.5){\line(-1,3){0.5}}
\put(-1,-2.5){\line(-1,5){0.5}}
\put(0,-0.5){\line(-1,1){0.5}}
\put(-1,-1.5){\circle*{0.133333}}

\put(0.5,3.5){\line(1,0){1}}
\put(1,3){\line(0,-1){0.5}}
\put(1,2.5){\line(1,2){0.5}}
\put(1,3){\line(1,1){0.5}}
\put(1,3){\line(-1,1){0.5}}
\put(0,4){\line(1,-1){0.5}}
\put(0,4.5){\line(1,-2){1}}
\put(0,4.5){\line(0,-1){0.5}}
\put(1,3){\circle*{0.133333}}
\put(1.5,3.5){\circle*{0.133333}}
\put(0,4){\circle*{0.133333}}
\put(0,4.5){\circle*{0.133333}}

\put(-0.5,3.5){\line(-1,0){1}}
\put(-1,3){\line(0,-1){0.5}}
\put(-1,2.5){\line(-1,2){0.5}}
\put(-1,3){\line(-1,1){0.5}}
\put(-1,3){\line(1,1){0.5}}
\put(0,4){\line(-1,-1){0.5}}
\put(0,4.5){\line(-1,-2){1}}
\put(0,4.5){\line(0,-1){0.5}}
\put(-1,3){\circle*{0.133333}}
\put(-1.5,3.5){\circle*{0.133333}}
\put(0,4){\circle*{0.133333}}
\put(0,4.5){\circle*{0.133333}}

\put(0.5,-3.5){\line(1,0){1}}
\put(1,-3){\line(0,1){0.5}}
\put(1,-2.5){\line(1,-2){0.5}}
\put(1,-3){\line(1,-1){0.5}}
\put(1,-3){\line(-1,-1){0.5}}
\put(0,-4){\line(1,1){0.5}}
\put(0,-4.5){\line(1,2){1}}
\put(0,-4.5){\line(0,1){0.5}}
\put(1,-3){\circle*{0.133333}}
\put(1.5,-3.5){\circle*{0.133333}}
\put(0,-4){\circle*{0.133333}}
\put(0,-4.5){\circle*{0.133333}}

\put(-0.5,-3.5){\line(-1,0){1}}
\put(-1,-3){\line(0,1){0.5}}
\put(-1,-2.5){\line(-1,-2){0.5}}
\put(-1,-3){\line(-1,-1){0.5}}
\put(-1,-3){\line(1,-1){0.5}}
\put(0,-4){\line(-1,1){0.5}}
\put(0,-4.5){\line(-1,2){1}}
\put(0,-4.5){\line(0,1){0.5}}
\put(-1,-3){\circle*{0.133333}}
\put(-1.5,-3.5){\circle*{0.133333}}
\put(0,-4){\circle*{0.133333}}
\put(0,-4.5){\circle*{0.133333}}

\end{picture}
\caption{The structure formed from a diamond-holder.} \label{structurefromholder}
\end{figure}
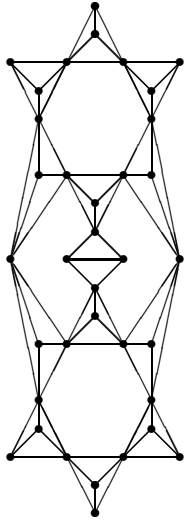

Hence,
any $C_{5}$ in $G^{*}$ must either be contained entirely within the structure formed from a diamond-holder
(see Figure~\ref{structurefromholder})
or must be contained within one snowflake.
However,
it can be seen that Figures~\ref{snowflake} and~\ref{structurefromholder} are both $C_{5}$-free,
and hence so is $G^{*}$.
\phantom{qwerty}
\setlength{\unitlength}{0.25cm}
\begin{picture}(1,1)
\put(0,0){\line(1,0){1}}
\put(0,0){\line(0,1){1}}
\put(1,1){\line(-1,0){1}}
\put(1,1){\line(0,-1){1}}
\end{picture} \\

As an aside,
let us remark that it is also possible to produce graphs attaining equality in Theorem~\ref{mainthm}
via alternative constructions to that used in Theorem~\ref{endthm}
(for example,
it is possible to obtain such graphs when $n \in \{14, 39, 54\}$).

\section{Concluding remarks} \label{concluding}

In this paper,
we have made a start on the topic of `extremal' planar graphs,
defining $\ex(n,H)$ to be the maximum number of edges possible in a planar graph on $n$ vertices 
that does not contain a given graph $H$ as a subgraph.
We observed
that the case when $H$ is a complete graph is straightforward,
and consequently focused on the more exciting examples of $C_{4}$ and $C_{5}$,
obtaining the tight bounds
$\ex(n,C_{4}) \leq \frac{15}{7}(n-2)$ for all $n \geq 4$ (in Section~\ref{C4})
and $\ex(n,C_{5}) \leq \frac{12n-33}{5}$ for all $n \geq 11$ (in Section~\ref{full}).

The next step would seem to be to look at larger cycles.
As the $C_{5}$ proof is rather intricate,
it may be impractical to try to obtain exact results,
but it would perhaps be reasonable to investigate the leading term in the extremal numbers 
and to work towards a general formula for this.

Of course,
there are also many other options to choose for the forbidden subgraph ---
the simplest would perhaps be $K_{4}$ minus one edge.
In particular,
it would be interesting to discover whether or not the chromatic number ever plays a role,
as in the Erdos-Stone Theorem.

Another angle to take on the problem would be to examine
which graphs have extremal number $3n-6$.
We have observed that $K_{4}$ certainly falls into this category,
and it would be interesting to develop necessary and sufficient conditions for all such graphs.

There are also various ways to extend the topic further.
For example,
one natural idea would be to investigate extremal problems on other surfaces,
such as the torus,
while another option would be to forbid several subgraphs at once.

There seems to be a large amount of uncharted territory here,
and it is hoped that this paper has made a useful start.

\section*{Acknowledgements}

This work was largely produced at the Ecole Polytechnique 
as part of the European Research Council funded project
ERC StG 208471 ExploreMaps,
and I am grateful for this support.
I would also like to thank the referees for some very helpful comments.

\end{document}